\documentclass[smallextended]{svpb}
\usepackage{amsmath, amssymb, theorem, latexsym,graphicx,fullpage}
\usepackage{amsmath, amssymb, theorem, latexsym,graphicx}
\usepackage{color}
\usepackage[pdftex]{hyperref}
\spnewtheorem{properties}{Properties}{\bf}{\it}  

\numberwithin{equation}{section} \numberwithin{figure}{section}

\renewcommand\theenumi{\roman{enumi}}

\allowdisplaybreaks

\def\cH {{\cal H}}
\def\cL {{\cal L}} 
\def\cN {{\cal N}} 
\def\cM {{\cal M}} 
\def\cQ {{\cal Q}} 

\def \N {{\mathbb N}}
\newcommand{\Sb}{\mathbb{S}}

\def \R {{\mathbb R}}

\newcommand{\sgn}[1]{\mathrm{sgn}({#1})}

\newcommand{\noi}{\noindent }
\newcommand{\noib}{\noindent $\bullet $~}
\newcommand{\ie}{i.e.}

\newcommand{\resp}{\emph{resp. }}
\newcommand{\etc}{\emph{etc\,}}

\title{A.~Stern's analysis of the nodal sets of some families of spherical harmonics revisited}

\titlerunning{Nodal sets of spherical harmonics}

\author{P. B\'erard \and B. Helffer}
\institute{P. B\'{e}rard \at Institut Fourier, Universit\'{e} de Grenoble and CNRS, B.P.74,
F 38402 Saint Martin d'H\`{e}res Cedex, France.\\\email{pierrehberard@gmail.com} \and B. Helffer \at Laboratoire de Math\'ematiques, Univ. Paris-Sud 11 and CNRS,
F 91405 Orsay Cedex, France, and
Laboratoire de Math\'{e}matiques Jean Leray, Universit\'{e} de Nantes.\\\email{Bernard.Helffer@math.u-psud.fr}}

\date{June 18, 2015}

\begin{document}
\maketitle

\begin{abstract}
In this paper, we revisit the analyses of Antonie Stern (1925)  and
Hans Lewy (1977) devoted to the construction of spherical harmonics
with two or three nodal domains. Our method yields sharp quantitative
results and a better understanding of the occurrence of bifurcations in the families of nodal sets.This paper is a natural continuation of
our critical reading of A. Stern's results for Dirichlet eigenfunctions in the square, see arXiv:14026054.
\end{abstract}


\section{Introduction}\label{S-intro}

Let $D$ be a regular bounded domain in $\R^n$. Let $\Delta$ be the
non-positive Laplacian with Dirichlet or Neumann boundary
conditions. We arrange the eigenvalues  $(\lambda_j)_{j\in \mathbb
N^*}$ of $- \Delta$ in increasing order,
$$\lambda_1 < \lambda_2 \le \lambda_3 \le \cdots .$$

Courant's 1923 celebrated nodal domain theorem \cite{Cou}, \cite[p.
452]{CH} states that an eigenfunction associated with the $n$-th
eigenvalue $\lambda_n$, has at most $n$ nodal domains. On the other
hand, an eigenfunction associated with $\lambda_n$, has at least two nodal domains when $n\ge 2$. The question remained of an eventual
lower bound for the number of nodal domains of an $n$-th
eigenfunction, as in the Sturm-Liouville theory. Antonie Stern's 1924 thesis \cite{St}, written under the supervision of Richard Courant, contains the following three results which provide a negative answer to this question.

\begin{theorem}\label{I-SP1}
Let $D$ be the unit square in $\R^2$, and $\Delta$ the non-positive
Laplacian with Dirichlet boundary conditions. Then, for any integer
$m$, there exists an eigenfunction $u$ of $-\Delta$, associated with
the eigenvalue $(4m^2+1)\pi^2$, whose nodal set inside the square
consists of a single simple closed curve. As a consequence, $u$ has
exactly two nodal domains.
\end{theorem}

\begin{theorem}\label{I-SP2odd} Let $\mathbb S^2$ be the unit
sphere in $\R^3$, and $\Delta$ the non-positive spherical Laplacian.
For any odd integer $\ell$, there exists a spherical
    harmonic of degree $\ell$, whose nodal set consists of a single
    simple closed curve. As a consequence, $u$ has exactly two nodal
    domains.
\end{theorem}

\begin{theorem}\label{I-SP2even}
Let $\mathbb S^2$ be the unit sphere in $\R^3$, and $\Delta$ the
non-positive spherical Laplacian.
 For any even integer $\ell \ge 2$, there exists a
    spherical harmonic of degree $\ell$, whose nodal set consists
    of two disjoint simple closed curves. As a consequence, $u$ has
    exactly three nodal domains.
\end{theorem}

Recall that the eigenfunctions of the spherical Laplacian are the spherical harmonics, \ie\, the restriction to the sphere $\mathbb S^2$ of the harmonic homogeneous polynomials in $\R^3$.\medskip

Theorem~\ref{I-SP1} is stated without proof in \cite[p. 455]{CH},
with a reference to Stern's thesis \cite{St}, and illustrated by two
figures taken from \cite{St}. Theorems ~\ref{I-SP2odd} and
\ref{I-SP2even} are apparently not mentioned in \cite{CH}. They
were rediscovered in 1977 by Hans Lewy \cite[Theorems~1 and 2]{Lew},
without any reference to A.~Stern. In the introduction of his paper,
Lewy also explains why a spherical harmonic of positive even degree has at least three nodal domains.

In \cite{St1}, we provide extracts from Stern's thesis, with
annotations and highlighting of the main assertions and ideas.
Stern's thesis is rather discursive. The main results are not stated
in theorem form, as above. They appear in the course of the
thesis, for example in the following citations (see Appendix~\ref{S-transl} for a translation into English)
 \cite[tags E1, K1, K2]{St1}:

\begin{quote}
[E1]~\ldots es l\"a\ss t sich beispielweise leicht zeigen, da{\ss} auf
der Kugel bei jedem Eigenwert die Gebietszahlen $2$ oder $3$
auftreten, und da{\ss} bei Ordnung nach wachsenden Eigenwerten auch beim
Quadrat die Gebietszahl $2$ immer wieder vor\-kom\-mt.\medskip

[K1]~Zun\"achst wollen wir zeigen, da{\ss} es zu jedem Eigenwert
Eigenfunktionen gibt, deren Nullinien die Kugelfl\"ache nur in zwei
oder drei Gebiete teilen. \ldots \ldots Die Gebietszahl zwei tritt somit
bei allen Eigenwerten  $\lambda_n = (2r+1)(2r+2)\, ~r=1,2,\cdots $
auf ; \\[3pt]
[K2]~ebenso wollen wir jetzt zeigen, da{\ss} die Gebietszahl
drei bei allen Eigenwerten $ \lambda_n =2r  (2r+1), ~r=1,2,\cdots$
immer wieder vorkommt.
\end{quote}

Stern's proofs are far from being complete, but she provides nice
geometric arguments \cite[tags I1-I3]{St1}, and figures.

\begin{quote}
[I1]~Legen wir die beiden Knotenliniensysteme \"ubereinander und
schraffieren wir die Gebiete, in denen beide Funktionen gleiches
Verzeichen haben, so kann die Knotenlinie der Kugelfunktion
$$
P^{2r+1}_{2r+1} (\cos \vartheta)\cos (2r+1) \varphi  + \mu
P_{2r+1}(\cos \vartheta)\,,\quad \mu >0
$$
nur in der nichtschraffierten Gebieten verlaufen\\[3pt]
[I3]~ und zwar f\"ur hinreichend kleine $\mu$ in beliebiger
Nachbarschaft der Knotenlinien von $$P^{2r+1}_{2r+1} (\cos
\vartheta)\cos (2r+1) \varphi,$$ d. h. der $2r+1$ Meridiane, da sich
bei stetiger \"Anderung von $\mu$ das Knotenliniensystem stetig
\"andert \ldots. \\[3pt]
[I2]~Da die Knotenlinie ferner durch die $2(2r+1)^2$ Schnittpunkte
der Nullinien der beiden obenstehenden Kugelfunktionen gehen mu{\ss}
\ldots
\end{quote}

\textbf{Sketch of Stern's proofs, and comparison with Lewy's paper \cite{Lew}}. A.~Stern starts from an eigenfunction $u$, whose nodal set can be completely described.  In the case of the  unit  square, this function $u$ is chosen to be $\sin(2m\pi x) \, \sin(\pi y) + \sin(\pi x) \,
\sin(2m\pi y)$. In the case of the sphere, it is chosen to be the
restriction to the sphere of the homogeneous harmonic polynomial
$W(x,y,z)=\Im(x+iy)^{\ell}$. A.~Stern then perturbs the eigenfunction $u$ by some eigenfunction $v$ (in the same eigenspace), looking at the
family $u^{\mu} = u+\mu v$ for $\mu$ small. The function $v$ is
chosen to be $\sin(\pi x) \, \sin(2m\pi y)$ in the case of the square, and a spherical harmonic whose nodal sets mainly contains latitude circles
(parallels) in the case of the sphere.

The main observation made by Stern is that for $\mu > 0$, the
nodal set $N(u^{\mu})$ satisfies
$$\cN \subset N(u^{\mu}) \subset \cN \cup \{u \, v < 0\},$$
where $\cN = N(u)\cap N(v)$ is the set of zeros common to $u$ and
$v$. In the case of the square, the connected components of the set
$\{u \, v \neq 0\}$ are small squares, whose vertices belong to
$\cN$. In the case of the sphere, they are square-like domains with
vertices in $\cN$, and triangle-like domains one of whose vertices
is the north or south pole, and the others belong to $\cN$. In both
cases, the domains form a kind of  grey/white checkerboard
(the connected open sets on which $u\, v$ is positive/negative) on
the unit square or on the sphere. The above inclusions say that the
nodal set $N(u^{\mu})$ contains $\cN$, and has to avoid the grey
squares, see \cite[tags I1, I2]{St1}. A.~Stern concludes, without proof, by saying that the nodal set $N(u^{\mu})$ deforms continuously, and remains close to the nodal set $N(u)$ when $\mu$ is small, \cite[tag
I3]{St1}.\medskip

To prove Theorems~\ref{I-SP2odd} and \ref{I-SP2even}, Lewy \cite{Lew} analytically determines how the nodal sets deform under small perturbations, first locally, and then globally. Continuity arguments were later explored in \cite{LeyD}.\smallskip

For a given degree $\ell$, A.~Stern \cite{St} and H.~Lewy
\cite{Lew} give examples of a spherical harmonic $h_o$ whose nodal
set is a simple closed regular curve, when the degree is odd; and  of a
spherical harmonic $h_e$ whose nodal set consists of two disjoint
simple closed regular curves, when the degree is even.

\emph{When the degree $\ell$ is odd}, Stern and Lewy  start from the
spherical harmonic $W$, whose expression in spherical coordinates is
given by $w(\vartheta,\varphi) = \sin^{\ell}(\vartheta)\, \sin(\ell
\varphi)$, and consider the family $W + \mu F$, where $F$ is a
spherical harmonic of degree $\ell$. Lewy \cite[Theorem 1]{Lew} only
requires that $F(p_+) > 0$, where $p_{+}$ is the north pole. Stern (see the proof of
Proposition~\ref{Ex1-P}) chooses $F$ to be the zonal spherical
harmonic given by $P_{\ell}(\cos \vartheta)$ in spherical
coordinates.  As we shall see, for $\mu$ small enough, $0$ is a
regular value of $W + \mu F$, and the nodal set $N(W+\mu F)$ is
connected. As a consequence, one has two-parameter families of
spherical harmonics whose nodal sets are simple closed regular
curves, and hence with two nodal domains.

\emph{When the degree $\ell \ge 2$ is even}, Stern and Lewy use
different functions. Stern starts from the spherical harmonic $W$, and perturbs it by the
spherical harmonic $V_{\alpha}$ whose expression in spherical
coordinates is given by $v_{\alpha}(\vartheta,\varphi) =
P_{\ell}^{1}(\cos\vartheta)\, \sin(\varphi - \alpha)$. As we shall
see, for $\mu$ small enough,  $0$ is a regular value of $W + \mu
V_{\alpha}$ (see the proof of Proposition~\ref{Ex2-P}), and the
nodal set $N(W+\mu V_{\alpha})$ has two connected components. This
construction gives us a three-parameter family of spherical
harmonics of even degree $\ell$, admitting $0$ as a regular value,
and whose nodal sets have two connected components, and
hence with three nodal domains. Lewy \cite[Theorem~2]{Lew} first constructs a spherical harmonic $F$
of the form $F(x,y,z) = x\, y\, F_1(x,y,z)$, where the function $F_1$ only depends on the distance to the north pole $p_{+}$. The
nodal set $N(F)$ consists of two orthogonal great circles through
the poles $p_{\pm}$, and $(\ell - 2)$ latitude circles. The nodal set
$N(F)$ has $4(\ell - 2) + 2$ singular points which are double
crossings. Lewy then constructs a spherical harmonic $G$, of degree
$\ell$, with ad hoc signs at the singular points of $N(F)$, so that
$F+\mu G$ desingularizes $N(F)$ when $\mu$ is small enough. He then
shows that, for $\mu$ small enough, $0$ is a regular value of $F+\mu
G$, and that the nodal set $N(F+\mu G)$ has two connected
components. This construction yields a two-parameter family of such
spherical harmonics.\medskip

As a matter of fact, one can prove that the set of spherical
harmonics of odd degree $\ell$ (\resp of even degree $\ell$), which
admit $0$ as a regular value, and whose nodal set is connected
(\resp whose nodal set has two connected components), is an open set
in the vector space $\cH_{\ell}$ of spherical harmonics of degree $\ell$.
This is a consequence of the inverse function
theorem and of the existence of the above examples.
\medskip

A.~Stern's proofs lack important details, while Lewy's paper is written in a rather condensed analytical style.  In \cite{BeHe}, we gave a complete geometric proof of Theorem~\ref{I-SP1}. In this paper, we give a complete geometric proof of quantitative versions of Theorems~\ref{I-SP2odd} and ~\ref{I-SP2even} (see  Propositions~\ref{Ex1-P} and \ref{Ex2-P} respectively). In each case, we start from Stern's geometric ideas, and we make a precise analysis of the possible local nodal patterns in the square-like and triangle-like domains. More precisely, we supplement Stern's ideas with,\vspace{-3mm}
\begin{enumerate}
    \item an analysis of the \emph{critical zeros} (a critical zero of a function $f$ is a point at which both $f$ and its differential $df$ vanish), showing in particular that $u^{\mu}$
    does not have any critical zero when $\mu \neq 0$ is small enough
    (Lemmas~\ref{Ex1-L1} and \ref{Ex2-L1});
    \item separation lemmas (Lemmas~\ref{Ex1-L2} and
    \ref{Ex2-L2}) to exclude certain local nodal patterns (in our proof, they replace Lewy's continuity arguments);
    \item an energy argument to show that a connected component of
    the set  $\{u\,v \neq 0\}$ cannot contain a simple closed nodal curve of $u^{\mu}$, see Properties~\ref{Ex1-P1}~(iii) and
    Properties~\ref{Ex2-P1}~(vii);
    \item classical properties of nodal sets of
    eigenfunctions, as summarized in \cite[Section~5.2]{BeHe}.
\end{enumerate}

Our proofs yield sharp quantitative results, Propositions~\ref{Ex1-P} and \ref{Ex2-P}, and a better understanding of the occurrence of bifurcations in the families of nodal sets, Lemmas~\ref{Ex1-L1} and \ref{Ex2-L1}. We in particular show that there exists some positive $\mu_c$ such that, for $0 < \mu < \mu_c$, the nodal set of $u^{\mu}$ is a regular $1$-dimensional submanifold of the sphere, while the nodal set of $u^{\mu_c}$ has self-intersections. \medskip

As far as Courant's theorem is concerned, another natural question is whether Courant's upper bound is sharp. For $2$-dimensional Euclidean domains with Dirichlet boundary condition, using the Faber-Krahn inequality in an essential way, {\AA}ke Pleijel \cite{Pl} proved that the number of nodal domains of an $n$-th
eigenfunction is asymptotically  less than $0.7\, n~$ (more precisely, less than $\gamma(2)\, n$, with $\gamma(2):=4\pi/\lambda (Disk_1)$, where $\lambda(Disk_1)$ is the lowest eigenvalue of the Dirichlet Laplacian in the disk of area $1$). As a corollary, one can conclude that Courant's theorem is sharp for finitely many eigenvalues only. Pleijel's result was later generalized to any compact Riemannian manifold $M$ (with Dirichlet boundary condition if $\partial M \neq \emptyset$), with a universal constant $\gamma(n) < 1$ depending only on the dimension $n$ of $M$, see \cite{Pe,BeMe}. The case of the Neumann boundary condition was also considered by Pleijel \cite{Pl} for the square, and
recently revisited in \cite{Pol} (in a more general setting: dimension 2, piecewise real analytic boundary), and in \cite{HP}. Starting in 2009, there has been a renewed interest for Courant's theorem in the context of
minimal partitions, and the investigation of the cases in which
Courant's theorem is sharp \cite{HHOT,HHOT1}, see Section~\ref{s5}.
These developments motivated \cite{BeHe} and motivate the present
paper.\medskip

The paper is organized as follows. In Section~\ref{S-Pre} we recall
some properties of spherical harmonics and Legendre functions. In
Sections~\ref{S-Ex1} and \ref{S-Ex2}, we give detailed geometric
proofs of Stern's second and third theorems for the sphere, with
quantitative statements (Propositions~\ref{Ex1-P} and \ref{Ex2-P}). In Section \ref{s5}, we recall the state of the art on the question of Courant sharpness for the sphere. In Appendix~\ref{S-Maple},
we provide some numerical computations of nodal sets of spherical
harmonics. In Appendix~\ref{S-transl}, we provide a translation into English of the citations from Stern's thesis.

\section{Preliminaries}\label{S-Pre}

\subsection{Spherical harmonics}\label{SS-SH}
Denote by $\mathbb S^2 := \{(x,y,z) \in \R^3 ~|~ x^2+y^2+z^2 = 1\}$
the round $2$-sphere. Given an integer $\ell \in \N$, we call
$\cH_{\ell}$ the vector space of spherical harmonics of degree
$\ell$ \ie, the restriction to the sphere of the harmonic
homogeneous polynomials in $3$ variables in $\R^3$. This is the
eigenspace of $-\Delta$ on $\mathbb S^2 $, associated with the
eigenvalue $\ell (\ell + 1)$. It has dimension $2\ell + 1$. Given a
spherical harmonic $h(\xi,\eta,\zeta)$ of degree $\ell$, with
$(\xi,\eta,\zeta) \in \Sb^2$, one can recover the harmonic
homogeneous polynomial $H$ it comes from as follows. Let $r =
(x^2+y^2+z^2)^{1/2}$. Then,
\begin{equation}\label{SH-1}
H(x,y,z) = r^{\ell} \, h(xr^{-1},yr^{-1},zr^{-1})\,.
\end{equation}
For simplicity, we shall henceforth identify the spherical harmonic
$h$ and the polynomial $H$. \medskip

The space $\cH_0$ is $1$-dimensional, associated with the eigenvalue
$0$. The space $\cH_1$ has dimension $3$. It is associated with the
eigenvalue $2$, and is generated by the coordinate functions $x, y$
and $z$ which have exactly two nodal domains. The space $\cH_2$ has
dimension $5$. It is associated with the eigenvalue $6$, and is
generated by the polynomials $yz, xz, xy, 2z^2-x^2-y^2$ and
$x^2-z^2$. It is easy to check that for $\mu > 0$, small enough, the
spherical harmonic $xy + \mu\, (2z^2-x^2-y^2)$ has exactly three
nodal domains: the nodal set of the spherical harmonic consists of
two simple closed curves given by the intersection of the sphere
with a right cylinder over a hyperbola in the $\{x,y\}$-plane.
Following A.~Stern, we shall later on consider a perturbation of the
degree $\ell$ spherical harmonic $\Im(x+iy)^{\ell}$.\medskip

We denote the north and south poles of $\mathbb S^2 $ respectively by $p_{+}
= (0,0,1)$ and $p_{-} = (0,0,-1)$.\medskip

By abuse of language, we shall call \emph{spherical coordinates} on
the sphere $\mathbb S^2$, the map
\begin{equation}\label{SH-4}
\left\{
\begin{aligned}
& E : [0,\pi]\times\R \to \Sb^2\,,\\
& E(\vartheta,\varphi)= (\sin\vartheta \, \cos\varphi, \sin\vartheta
\, \sin\varphi, \cos\vartheta)\,,
\end{aligned}
\right.
\end{equation}
where $\vartheta$ is the \emph{co-latitude}, and $\varphi$ the
\emph{longitude}.\\
 The map $E$ is a diffeomorphism from
$(0,\pi)\times(\varphi_0, 2\pi+\varphi_0)$ onto $\Sb^2\setminus
M_{\varphi_0}$, where $M_{\varphi_0}= E([0,\pi],\varphi_0)$ is the
meridian from $p_{+}$ to $p_{-}$ with longitude $\varphi_0$. To
cover $\Sb^2\setminus\{p_{\pm}\}$, we will work in
$(0,\pi)\times\R_{2\pi}$ \ie, modulo $2\pi$ in the $\varphi$
variable ($\mathbb R_{2\pi} = \mathbb R/ (2\pi \mathbb Z)$).

The map $E$ can be viewed as the polar coordinates in the
exponential map $exp_{p_{+}}$, which sends the disk $D(0,\pi)$, with
center $0$ and radius $\pi$ in $T_{p_{+}}\Sb^2$ (the tangent plane to the sphere at $p_{+}$), onto
$\Sb^2\setminus \{p_{-}\}$ diffeomorphically. The variable
$\vartheta$ is the distance to the north pole.

In the spherical coordinates, the antipodal map is given by
$(\vartheta,\varphi) \to (\pi-\vartheta,\pi+\varphi)$.

In the sequel, we will illustrate the proofs by figures representing
the nodal patterns viewed through the exponential map \ie, in the
disk $D(0,\pi)$. In the figures, the outer circle always represents the circle of radius $\pi$ \ie, the cut-locus of the north pole.

Using the spherical coordinates, the spherical harmonics can be
described in terms of Legendre functions and polynomials. In the
next section, we fix some notation, and recall useful properties
of Legendre functions and polynomials.

\subsection{Legendre functions and polynomials}\label{SS-P}

The $(2\ell + 1)$-dimensional vector space $\cH_{\ell}$ of spherical
harmonics of degree $\ell$ admits the basis,
\begin{equation}\label{P-1}
P_{\ell}(\cos\vartheta),
~P_{\ell}^{m}(\cos\vartheta)\,\cos(m\varphi),
~P_{\ell}^{m}(\cos\vartheta)\,\sin(m\varphi)\,,
\end{equation}
where $m$ is an integer $1 \le m \le \ell$, $P_{\ell}$ the Legendre
polynomial of degree $\ell$, and $P_{\ell}^{m}$ the Legendre
functions. We use the notation and the normalizations of \cite{MOS}.

\begin{properties}\label{P-P1}
For $0 \le m \le \ell$, the Legendre function $P_{\ell}^{m}$
satisfies the differential equation,
\begin{equation}\label{P-2}
\left( (1-t^2) P'(t)\right)' + \left( \ell (\ell + 1) -
\frac{m^2}{1-t^2}\right) P(t) = 0\,.
\end{equation}
When $m$ is $0$, $P_{\ell}^m = P_{\ell}$\,. Furthermore, the
following properties hold. \vspace{-3mm}
\begin{enumerate}
\item Identities.
\begin{equation}\label{P-3a}
P_{\ell}(t)  = \frac{1}{2^{\ell} \ell!}
\left(\frac{d}{dt}\right)^{\ell}(t^2-1)^{\ell}\,,
\end{equation}
\begin{equation}\label{P-3b}
P_{\ell}^m(t) =
(-1)^m(1-t^2)^{m/2}\left(\frac{d}{dt}\right)^{m}P_{\ell}(t)\,,
\end{equation}
and
\begin{equation}\label{P-3c}
(1-t^2)\,P_{\ell}'(t)  = \ell \, P_{\ell -1}(t) - \ell\,t\,
P_{\ell}(t)\,.
\end{equation}
In particular, $P_{\ell}^{\ell}(\cos\vartheta) = C_{\ell}
\sin^{\ell}(\vartheta)\,$, where $C_{\ell}$ is a constant.

\item The polynomial $P_{\ell}$ has degree $\ell$, the same
parity as the integer $\ell$, and satisfies
\begin{equation}\label{P-pm1}
P_{\ell}(1) = 1\,, \, P_{\ell}(-1) = (-1)^{\ell}\,,
\end{equation}
and
\begin{equation}\label{P-sup}
\sup_{[-1,1]}|P_{\ell}(t)|=1\,.
\end{equation}

\item The polynomial $P_{\ell}(t)$ has $\ell$ simple roots $t_j(\ell)$ ($j=1, \cdots,\ell$) in the interval $(-1,1)$, enumerated in
decreasing order. We write these roots as $t_j(\ell)=
\cos(\vartheta_j(\ell))$, with
\begin{equation}\label{P-r}
0 < \vartheta_1(\ell) < \vartheta_2(\ell) < \cdots < \vartheta_{\ell-1} (\ell) <
\vartheta_{\ell} (\ell) < \pi \,.
\end{equation}
The derivative $P'_{\ell}(t)$ has $(\ell - 1)$ simple roots which we
write as $\cos\left(\vartheta'_j(\ell)\right)$. They satisfy
\begin{equation}\label{P-r1}
0 < \vartheta_1 (\ell) < \vartheta'_1(\ell)  < \vartheta_2 (\ell) < \cdots <
\vartheta_{\ell-1}(\ell)  < \vartheta'_{\ell-1}(\ell)  < \vartheta_{\ell}(\ell) < \pi
\,.
\end{equation}
Note that the values $\vartheta_j(\ell)$ and $\vartheta'_j(\ell) $ are
symmetrical with respect to $\frac{\pi}{2}$. As a consequence, for
$\ell$ odd, $P_{\ell}(0) = 0$ and $P'_{\ell}(0) \neq 0$, and for
$\ell$ even, $P_{\ell}(0) \neq 0$ and $P'_{\ell}(0) = 0$.

\item The polynomials $P_{\ell}$ and $P_{\ell - 1}$ have no common
zero. More precisely, the zeros of $P_{\ell}$ and $P_{\ell-1}$ are
intertwined:
\begin{equation}\label{P-z}
0 < \vartheta_{1}(\ell) < \vartheta_{1}(\ell-1) <
\vartheta_{2}(\ell) < \cdots < \vartheta_{\ell-1}(\ell) <
\vartheta_{\ell-1}(\ell-1) < \vartheta_{\ell}(\ell) < \pi \,.
\end{equation}
\item For the zeros of $P_{\ell}$, one has the inequalities,
\begin{equation}\label{P-Sz}
\frac{2 j -1}{2\ell +1}  \pi < \vartheta_j (\ell)  < \frac{2 j}{2\ell +1}\pi \,,
\, \mbox{ for } j=1,\cdots,\ell\,.
\end{equation}
\item Call $p_j(\ell), ~1 \le j \le [\frac{\ell}{2}]$, the local
maxima of $|P_{\ell}(t)|$, when $t$ decreases from $1$ to $0$. Then,
\begin{equation}\label{P-Szlm}
0 < p_{[\frac{\ell}{2}]}(\ell) < \cdots < p_2(\ell) < p_1(\ell) < 1\,.
\end{equation}
Here $[\cdot]$ denotes the integer part.
\item For the derivative of the Legendre polynomial $P_{\ell}(t)$,
one has the inequality,
\begin{equation}\label{P-Szmax}
|P'_{\ell}(t)| \le \frac{\ell (\ell + 1)}{2} \text{~~for~~} - 1 \le
t \le 1 \,,
\end{equation}
where the equality is achieved for $\ell = 0, 1$ and when $\ell \ge
2$, for $t = \pm 1$.
\end{enumerate}
\end{properties}

For these properties, we refer to \cite[Chapters~IV and V]{MOS} and
to \cite{Sz}. In particular, Properties (v)-(vii) can be found in
\cite{Sz}, \resp under Theorems~6.21.2, 7.3.1, and Inequality
(7.33.8). \medskip

\textbf{Remarks}.  (i) Using \eqref{P-3c}, one can prove that for $1
\le j \le \ell-1$,
\begin{equation}
\vartheta_j(\ell) < \vartheta_j(\ell-1) < \vartheta'_j(\ell)\,.
\end{equation}

(ii) One can  relate the asymptotic behavior of $\vartheta_1(\ell)$
as $\ell \rightarrow +\infty$, to the first zero $ j_{0,1} $ of the
zero-th Bessel function $J_0$, \cite[Theorem~8.1.2]{Sz},
\begin{equation} \label{asvartheta}
\vartheta_1(\ell) \sim { j_{0,1}} /\ell\,.
\end{equation}

(iii)  One can also relate the asymptotic behavior of $p_1(\ell)$ as
$\ell \rightarrow +\infty$, to the first zero $ j_{1,1}$ of the
Bessel function $J_1=-J'_0$,
\begin{equation}\label{asp1}
p_1(\ell) \sim - J_0( j_{1,1})\,.
\end{equation}


\section{Stern's first theorem: odd case}\label{S-Ex1}

The purpose of this section is to prove Theorem~\ref{I-SP2odd}. As a
matter of fact, we shall prove a more quantitative result,
Proposition~\ref{Ex1-P1}, which implies the theorem. We use Stern's
ideas sketched in the introduction.

Fix an integer $\ell \in \N$, without any parity assumption for the
time being. We work in the spherical coordinates \eqref{SH-4}, with $(\vartheta,\varphi) \in [0,\pi] \times \R_{2\pi}$.

\subsection{Notation}\label{SS-Ex1Not}

Up to scaling, there is a unique spherical harmonic $Z_{\ell}$, of
degree $\ell$, which is invariant under the rotations about the
$z$-axis. Viewed in the spherical coordinates, this zonal spherical
harmonic is given by $P_{\ell}(\cos \vartheta)$. Let
$\vartheta_1(\ell) < \vartheta_2(\ell) < \cdots <
\vartheta_{\ell}(\ell)$ be the zeros of the function $\vartheta \to
P_{\ell}(\cos \vartheta)$ in the interval $(0,\pi)$, see
Properties~\ref{P-P1}~(iii). The nodal set of the spherical harmonic
$Z_{\ell}$, denoted $N(Z_{\ell})$, consists of precisely $\ell$
latitude circles (parallels),
\begin{equation*}
L_{i} := \{(\vartheta, \varphi) ~|~ \vartheta =
\vartheta_i(\ell)\}\,, ~1 \le i \le \ell\,.
\end{equation*}
They determine sectors on the sphere,
\begin{equation*}
\cL_{i} := \{(\vartheta, \varphi) ~|~ \vartheta_{i}(\ell) <
\vartheta < \vartheta_{i+1}(\ell)\}\,, ~0 \le i \le \ell\,,
\end{equation*}
where $\vartheta_0(\ell) = 0$ and $\vartheta_{\ell + 1}(\ell) =
\pi$. In the sector $\cL_i$, the function $Z_{\ell}$ has the sign of
$(-1)^i\,$.
\smallskip

Call $W_{\ell}$ the spherical harmonic of degree $\ell$ obtained by
restricting the harmonic homogeneous polynomial $\Im(x+iy)^{\ell}$
to the sphere. Viewed in spherical coordinates, this spherical
harmonic is given by $\sin^{\ell}(\vartheta) \, \sin(\ell \varphi)$.
Its nodal set $N(W_{\ell})$ consists of $\ell$ great circles of
$\mathbb S^2 $ \ie, of $2\ell$ meridians,
\begin{equation*}
M_j = \{(\vartheta, \varphi) ~|~ \varphi = j\frac{\pi}{\ell}\}\,, ~
0 \le j \le 2\ell - 1\,.
\end{equation*}
They determine sectors on the sphere,
\begin{equation*}
\cM_j = \{(\vartheta, \varphi) ~|~ j \frac{\pi}{\ell} < \varphi <
(j+1)\frac{\pi}{\ell}\}\, ~0 \le j \le 2\ell - 1 \,.
\end{equation*}
In the sector $\cM_j$ the function $W_{\ell}$ has the sign of
$(-1)^j$.

Note that these meridians meet at the north and south poles $p_{\pm}$
which are the only singular points of the nodal set $N(W_{\ell})\,$.
\medskip

The intersection $$\cN = N(Z_{\ell}) \cap N(W_{\ell})$$ is the
finite set of zeros common to $Z_{\ell}$ and $W_{\ell}$.

We call $q_{i,j}$ the intersection point of $L_i$ with $M_j\,$, $1
\le i \le \ell\,$, $0 \le j \le 2\ell-1\,$, so that
$$
\cN = \{q_{i,j}, ~1 \le i \le \ell, ~0 \le j \le 2\ell-1\}\,.
$$

For $0 \le i \le \ell$ and $0 \le j \le 2\ell-1$, we introduce the
sets $\cQ_{i,j} = \cL_i \cap \cM_j$, which are the connected
components of the open set $\{Z_{\ell} \, W_{\ell} \not = 0\}$. In
$\cQ_{i,j}$ the sign of the function $Z_{\ell} \, W_{\ell}$ is
$(-1)^{i+j}\,$.

The sets $\cQ_{i,j}$ form a grid patterns over the sphere, and
following the idea of A.~Stern, they can be colored according to the
sign of the function $Z_{\ell} \, W_{\ell}$ thus forming a \emph{checkerboard}.

Finally, for $0 \le j \le 2\ell -1\,$, we introduce the meridian
$B_j$,
\begin{equation}\label{Ex1-mb}
B_j = \{(\vartheta, \varphi) ~|~ \varphi = (j+\frac{1}{2})
\frac{\pi}{\ell}\}\,,
\end{equation}
which bisects the sector $\cM_j$.

\begin{figure}[!ht]
\begin{center}
\includegraphics[width=0.7\linewidth]{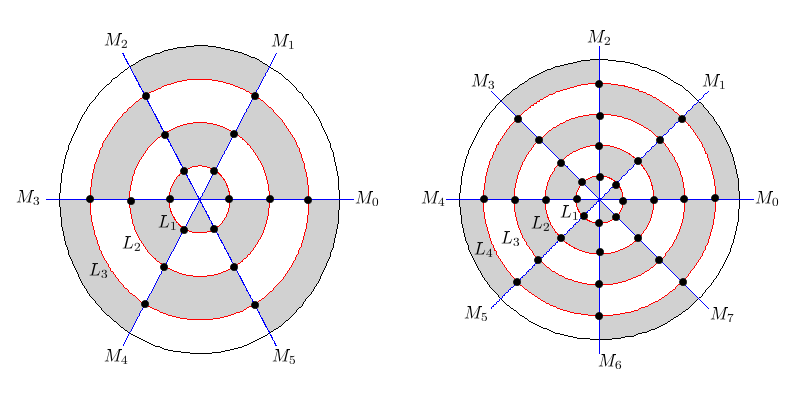}
\caption{Checkerboards in the cases $\ell = 3$, and $\ell = 4$}
\label{FEx1-nodal1}
\end{center}
\end{figure}

Figure~\ref{FEx1-nodal1} displays the latitude circles and the
meridians viewed through the exponential map at the north pole
$p_{+}$, in the cases $\ell = 3$ and $\ell = 4$. The common zeros of
$W_{\ell}$ and $Z_{\ell}$ are the big dots. The coloring white/grey
illustrates the sign of $Z_{\ell}\, W_{\ell} $. The
outer circle is mapped to the south pole by the exponential map.

\subsection{The family $H^{\mu,\ell}$}\label{SS-Ex1Fam}

Following Stern \cite{St}, we consider the one-parameter family of
spherical harmonics,
\begin{equation}\label{Ex1-1}
H^{\mu,\ell} = W_{\ell} + \mu \, Z_{\ell}\,,
\end{equation}
which may be written in spherical coordinates as
\begin{equation}\label{Ex1-2}
h^{\mu,\ell}(\vartheta,\varphi) =
\sin^{\ell}(\vartheta)\,\sin(\ell\varphi) + \mu \,
P_{\ell}(\cos\vartheta)\,.
\end{equation}

Note that
\begin{equation}\label{Ex1-sgnmu}
h^{-\mu,\ell}(\vartheta,\varphi) = - h^{\mu,\ell}(\vartheta,\varphi
+ \frac{\pi}{\ell})\,,
\end{equation}
for $(\vartheta,\varphi) \in (0,\pi)\times \R_{2\pi}$. It follows
that we can restrict to the case $\mu > 0$. We shall do so for the
remaining part of Section~\ref{S-Ex1}.

\subsubsection{Critical zeros of $H^{\mu,\ell}$}\label{SSS-Ex1-cz}

We call \emph{critical zero} of a function a point which is both a
zero and a critical point.

According to Properties~\ref{P-P1}~(ii), $H^{\mu,\ell}(p_{+}) =
\mu$, and $H^{\mu,\ell}(p_{-}) = (-1)^{\ell}\mu$, and hence the
north and south poles do not belong to the nodal set
$N(H^{\mu,\ell})$ when $\mu \neq 0$. As a consequence, for $\mu >
0$, the critical zeros of $H^{\mu,\ell}$ are located in
$\Sb^2\setminus \{p_{\pm}\}$, and we can look for them in the spherical
coordinates $(\vartheta,\varphi) \in \, (0,\pi)\times \R_{2\pi}$.

For $\mu > 0$, the point $(\vartheta,\varphi)$ corresponds to a
critical zero of $H^{\mu,\ell}$, if and only if
\begin{equation}\label{Ex1-3}
\begin{split}
& h^{\mu,\ell}(\vartheta,\varphi) = 0\,,\\
& \partial_{\vartheta}h^{\mu,\ell}(\vartheta,\varphi) = 0\,,\\
& \partial_{\varphi}h^{\mu,\ell}(\vartheta,\varphi) = 0\,.\\
\end{split}
\end{equation}
This is equivalent for $(\vartheta,\varphi)$ to satisfy the
relations,
\begin{equation}\label{Ex1-4a}
\begin{split}
& \cos(\ell \varphi) = 0 \text{~~\ie, ~} \sin(\ell \varphi) = \pm
1\,,\\
& \pm \, \sin^{\ell}(\vartheta) + \mu \, P_{\ell}(\cos\vartheta) = 0\,,\\
& \pm \, \ell\, \cos\vartheta \, \sin^{\ell-1}(\vartheta) - \mu\,
\sin \vartheta \, P_{\ell}'(\cos\vartheta) = 0 \,.
\end{split}
\end{equation}
Plugging \eqref{P-3c} into the third line of the above system, it
follows that, for $\mu \neq 0$, \eqref{Ex1-4a} is equivalent to
\begin{equation}\label{Ex1-4b}
\begin{split}
& \cos(\ell \varphi) = 0 \text{~~\ie, ~} \sin(\ell \varphi) = \pm
1\,,\\
& \frac{1}{\mu} = \mp \frac{P_{\ell}(\cos \vartheta)}{\sin^{\ell}(\vartheta)}\,,\\
& P_{\ell-1}(\cos \vartheta) = 0\,.
\end{split}
\end{equation}
%
%
By Properties~\ref{P-P1}~(iii)-(iv), the last equation in
\eqref{Ex1-4b} has exactly $(\ell-1)$ simple roots in $[0,\pi]$. We
denote them by, $\vartheta_1(\ell-1) < \ldots <
\vartheta_{\ell-1}(\ell-1)$. They are symmetrical with respect to
$\frac{\pi}{2}$ due to the parity of $P_{\ell-1}$.

It follows that, for $\mu > 0$, the only possible critical zeros of
the spherical harmonic $H^{\mu,\ell}$ are given in spherical
coordinates by the points $\left(
\vartheta_i(\ell-1),(j+\frac{1}{2})\frac{\pi}{\ell} \right)$ for $1
\le i \le \ell - 1$ and $0 \le j \le 2\ell-1$. These points can only
occur as critical zeros for finitely many values of $\mu$, given by
the second equation in \eqref{Ex1-4b}. Away from these values of
$\mu$, the spherical harmonic $H^{\mu,\ell}$ has no critical zero.

Since we restrict to $\mu > 0$, the \emph{critical} values of $\mu$
are given by
\begin{equation}\label{Ex1-ev}
\mu_i (\ell) =
\frac{\sin^{\ell}\left(\vartheta_i(\ell-1)\right)}{|P_{\ell}\left(\cos\vartheta_i(\ell-1)\right)|}\,,
\end{equation}
for $1 \le i \le \ell - 1$. \\
They are well-defined because the denominators do not vanish, since
the zeros of $P_{\ell}$ and $P_{\ell - 1}$ are intertwined, see
Properties~\ref{P-P1}~(iv). For the value $\mu_i(\ell)$, the
spherical harmonic $H^{\mu_i(\ell),\ell}$ has finitely many critical zeros
which are well determined by equations \eqref{Ex1-4b}. Note that the
values $\mu_i(\ell)$ are positive.

Taking the parity of the Legendre polynomials into account, it
suffices to consider the values $\mu_i(\ell)$ for $1 \le i \le
[\frac{\ell}{2}]$, where $[\frac{\ell}{2}]$ denotes the integer part
of $\frac{\ell}{2}$. We summarize the preceding discussion in
the following lemma.

\begin{lemma}\label{Ex1-L1}
Assume $\mu > 0$, and define $\mu_c (\ell)> 0$ to be the infimum
\begin{equation}\label{Ex1-mua}
\mu_c(\ell)  = \inf_{1 \le i \le [\frac{\ell}{2}]}\mu_i(\ell)\,,
\end{equation}
where the positive values $\mu_i(\ell)$ are given by \eqref{Ex1-ev}.

The spherical harmonic
$$
H^{\mu,\ell} = W_{\ell} + \mu \, Z_{\ell}
$$
does not vanish at the north and south poles. Except for the values
$\{\mu_i(\ell)\}_{i=1}^{[\frac{\ell}{2}]}$, $H^{\mu,\ell}$ has no critical
zero. In particular, for $0 < \mu < \mu_c(\ell)$, the function
$H^{\mu,\ell}$ has no critical zero, its nodal set is a
$1$-dimensional submanifold of the sphere, and hence consists of
finitely many disjoint regular simple closed curves.
\end{lemma}

The last assertion in the lemma follows from the fact that
self-intersections in the nodal set of an eigenfunction correspond
to critical zeros, see \cite{BeHe}, Section~5.2. \medskip

\textbf{Remark}. One can easily estimate $\mu_c(\ell) $ from below.
Using Properties~\ref{P-P1}~(ii) and (vi), one finds that
$$ \mu_i (\ell) > \frac{ \sin^{\ell}(\vartheta_i(\ell-1))}{p_i(\ell)}
> \frac{ \sin^{\ell}(\vartheta_1(\ell-1))}{p_1(\ell)}\,.$$
%

Hence,
\begin{equation}\label{lbmu}
 \mu_c(\ell)  \ge
\frac{\sin^{\ell}(\vartheta_1(\ell-1))}{p_1(\ell)}\,.
\end{equation}


Note that one has,

\begin{equation}\label{lbmu2}
\mu_c(\ell)  \ge  \min \left( \mu_1(\ell),
\frac{\sin^{\ell}(\vartheta_2(\ell))}{p_2(\ell)}\right)\,.
\end{equation}

 One can obtain the asymptotics of $\mu_1(\ell)$ as $\ell$
tends to infinity. Recall that $\frac{\pi}{2\ell+1} \le
\vartheta_1(\ell) \le \frac{2\pi}{2\ell+1}$, and use Hilb's formula
\cite[Theorems~8.21.6]{Sz},

$$
P_{\ell} (\cos \vartheta) = \left(\frac{\vartheta}{\sin
\vartheta}\right)^\frac{1}{2} \, J_0 \left((\ell+\frac{1}{2})
\vartheta\right) + R(\vartheta)\,,
$$
with
$$
R(\vartheta) = \mathcal O(\theta^2)\,,\, \mbox{ if } |\vartheta |
\leq C/{\ell}\,.
$$

It follows that
$$
\vartheta_1(\ell)= \frac{j_{0,1} }{\ell+ \frac{1}{2}} + \mathcal O
(1/{\ell}^3)\,,
$$
where $j_{0,1} $ is the least positive zero of the Bessel function
$J_0$.

Compute
\begin{equation*}
\begin{split}
P_{\ell}\left(\cos \vartheta_1(\ell-1)\right) &= \left(1 + \mathcal
O(\frac {1}{{\ell}^2}) \right)\, J_0 \left( (\ell+ \frac{1}{2})
\frac{j_{0,1} }{\ell-\frac{1}{2}}\right) + \mathcal O (\frac
{1}{\ell^2})\\
& = \frac{j_{0,1}  J'_0 (j_{0,1})}{(\ell -\frac{1}{2})}  +
\mathcal O(\frac{1}{\ell^2})\,.
\end{split}
\end{equation*}

Recall that
$$
\mu_1(\ell) = \frac{\sin^\ell \left(\vartheta_1(\ell-1)\right)}{| P_\ell
\left(\cos \vartheta_1(\ell-1)\right)|}\,.
$$

Observe that
$$
\sin\left(\vartheta_1(\ell-1)\right) =  \frac{ j_{0,1}}{\ell -
\frac{1}{2}} + \mathcal O (1/\ell^3)\,.
$$
Taking the power $\ell$, the remainder term does not change the main
term of the asymptotics,

$$
\sin^{\ell} \left(\vartheta_1(\ell-1)\right) \sim \left( \frac{
j_{0,1}}{\ell - \frac{1}{2}}\right)^\ell\,.
$$

Finally, we obtain

$$
\mu_1(\ell) \sim \left( \frac{j_{0,1} }{\ell - \frac{1}{2}}
\right)^{\ell-1}\frac{1}{ |J'_0({ j_{0,1}})|}\,.
$$

It turns out that the second term in the right hand side of
\eqref{lbmu2} is asymptotically bigger than the first one. It
follows that the preceding formula holds with $\mu_1(\ell)$ replaced
by $\mu_c(\ell)$ as well.

\subsubsection{A separation lemma for $N(H^{\mu,\ell})$}\label{SSS-Ex1-sl}

For $0 \le j \le 2\ell - 1\,$, we look at the function
$H^{\mu,\ell}$ restricted to the meridian $B_j\,$. Let
$$b^{\,\mu, \ell, j} = H^{\mu,\ell}|_{B_j}\,,$$
\ie,
$$b^{\,\mu, \ell, j}(\vartheta) = (-1)^{j}\, \sin^{\ell}(\vartheta)
+ \mu \, P_{\ell}(\cos \vartheta)\,.$$

Recall the notation  $\vartheta_1(\ell)$ and $p_1(\ell)$,
see Properties~\ref{P-P1}, (iii) and (vi).

\begin{lemma}\label{Ex1-L2}
The functions $b^{\,\mu, \ell, j}$ satisfy the following
properties.\vspace{-3mm}
\begin{enumerate}
\item  For $0 < \mu < \mu_c(\ell)$, the function
$b^{\,\mu, \ell, j}$ does not vanish in
$[\vartheta_1(\ell), \pi -\vartheta_1(\ell)]\,$.
\item When $\ell$ and $j$ are even, $b^{\,\mu, \ell, j}(\vartheta) >
0$ in $[0,\vartheta_1(\ell)] \cup [\pi-\vartheta_1(\ell),\pi]$.
\item When $\ell$ is even and $j$ odd, the function $b^{\,\mu, \ell,
j}(\vartheta)$ vanishes exactly once in each interval
$\left(0,\vartheta_1(\ell)\right)$ and
$(\pi-\vartheta_1(\ell),\pi)$.
\item When $\ell$ is odd and $j$ even, the function $b^{\,\mu, \ell,
j}(\vartheta)$ is positive in $[0,\vartheta_1(\ell)]$ and vanishes
exactly once in $(\pi-\vartheta_1(\ell),\pi)$.
\item When $\ell$ and $j$ are odd, the function $b^{\,\mu, \ell,
j}(\vartheta)$ vanishes exactly once in
$\left(0,\vartheta_1(\ell)\right)$ and is negative in
$[\pi-\vartheta_1(\ell),\pi]$.
\end{enumerate}\vspace{-3mm}
The previous items describe the possible intersections of the nodal
set $N(H^{\mu,\ell})$ with the meridian $B_j$ which bisects the
sector $\cM_j$.
\end{lemma}

\textbf{Proof}. (i)~ The function $b^{\,\mu, \ell, j}$
vanishes in the given interval if and only if
$$\frac{P_{\ell}(\cos\vartheta)}{\sin^{\ell}(\vartheta)} =
\frac{(-1)^{j+1}}{\mu}\,.$$ Call $\beta_1(\vartheta)$ the function in
the left hand side. Using \eqref{P-3c}, we obtain its derivative,
$\beta_1'(\vartheta) = - \, \frac{\ell \,
P_{\ell-1}(\cos\vartheta)}{\sin^{\ell+1}(\vartheta)}$. It follows
that the local extrema of the function $\beta_1$ are  achieved
at the values $\vartheta_i(\ell-1)$, and that
$$
b^{\,\mu, \ell, j}(\vartheta_i(\ell-1)) =
(-1)^j\,\sin^{\ell}(\vartheta_i(\ell-1))\,\Big[ 1 + (-1)^j \mu
\frac{P_{\ell}(\cos
\vartheta_i(\ell-1))}{\sin^{\ell}(\vartheta_i(\ell-1))}\Big]\,.
$$
Using \eqref{Ex1-ev}, the assertion follows. To prove assertions
(ii)-(v), we look at the signs of the functions
$P_{\ell}(\cos\vartheta)$ and $P'_{\ell}(\cos\vartheta)$ in the
intervals $\left(0,\vartheta_1(\ell)\right)$ and
$(\pi-\vartheta_1(\ell),\pi)$ in order to determine the signs of
$b^{\,\mu, \ell, j}(\vartheta)$ and $\partial_{\vartheta}b^{\,\mu,
\ell, j}(\vartheta)$. We leave the details to the reader. \hfill
$\square$ \medskip

\subsubsection{General properties of $N(H^{\mu,\ell})$}\label{SSS-Ex1-gp}

We now state simple general properties of the nodal set of the
spherical harmonic $H^{\mu,\ell} = W_{\ell} + \mu \, Z_{\ell}\,$. We
use the notation of Subsection~\ref{SS-Ex1Not}.

\begin{properties}\label{Ex1-P1}
For $\mu > 0$, the nodal sets of the spherical harmonics
$H^{\mu,\ell}$ share the following properties.\vspace{-3mm}
\begin{enumerate}
    \item The nodal set of $H^{\mu,\ell}$ satisfies
    $$\cN \subset N(H^{\mu,\ell}) \subset \cN \cup \{Z_{\ell} \, W_{\ell} <
    0\}.$$ This means that a point in the nodal set of $H^{\mu,\ell}$ is
    either one of the points in $\cN$, or a point in some open domain
    $\cQ_{i,j}$, with $(-1)^{i+j} = - 1$.
    \item The nodal set of $H^{\mu,\ell}$ near each point
    $q_{i,j} \in \cN$ consists of a single regular arc which is transversal to the
    latitude circle $L_i$ and to the meridian $M_j$. In other words, an arc in the
    nodal set inside some domain $\cQ_{i,j}$, with $(-1)^{i+j} = - 1$, can only exit
    $\cQ_{i,j}$ through a point in $\cN$ (a vertex), and cannot cross the boundary of
    $\cQ_{i,j}$ elsewhere.
    \item For $0 < \mu < \mu_c(\ell)$ (defined in Lemma~\ref{Ex1-L1}), no connected component
    of the nodal set $N(H^{\mu,\ell})$ can be entirely contained in some $\cQ_{i,j}\,$.
\end{enumerate}
\end{properties}

\textbf{Proof}. Property~(i) is clear. Property~(ii) follows from
the fact that both $\partial_{\vartheta}h^{\mu,\ell}$ and
$\partial_{\varphi}h^{\mu,\ell}$ do not vanish at the points
$(\vartheta_i(\ell),j\frac{\pi}{\ell})$. Indeed, $\mu > 0$, and $P_{\ell}$
and $P'_{\ell}$ have no common zero according to
Properties~\ref{P-P1}~(iv). For the proof of Property~(iii), we
observe that the assumption on  $\mu$ implies that the nodal set
$N(H^{\mu,\ell})$ is a $1$-dimensional submanifold of the sphere
\ie, a finite collection of disjoint regular simple closed curves.
Assume that one such closed curve is entirely contained in some
domain $\cQ_{i,j}$. This domain $\cQ_{i,j}$  would contain
some nodal domain of $H^{\mu,\ell}$, and hence its first Dirichlet
eigenvalue $\lambda$ would be strictly less than $\ell (\ell + 1)$,
the eigenvalue associated with $H^{\mu,\ell}$. On the other-hand,
since $\cQ_{i,j}$ is contained in one of the nodal domains of
$Z_{\ell}$ (or of $W_{\ell}$), its first eigenvalue satisfies
$\lambda > \ell (\ell + 1)\,$, a contradiction. \hfill
$\square$

\subsubsection{Local nodal patterns for $H^{\mu,\ell}$}\label{SSS-Ex1-lnp}

Assume that $0 < \mu < \mu_c(\ell)$.

\begin{figure}
\begin{minipage}[c]{.46\linewidth}
\includegraphics[width=6cm]{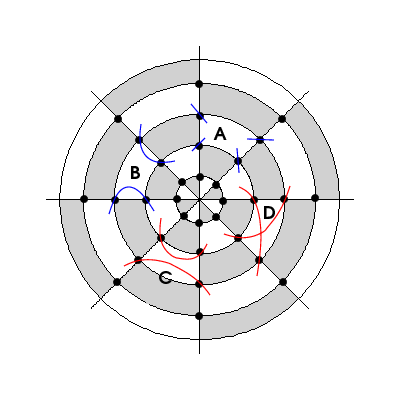}
\caption{Patterns away from the poles} \label{FEx1-np1}
\end{minipage} \hfill
\begin{minipage}[c]{.46\linewidth}
\includegraphics[width=6cm]{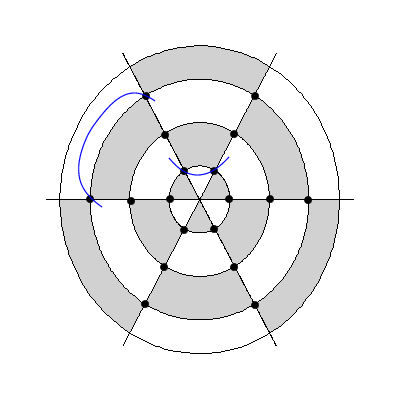}
\caption{Patterns near the poles} \label{FEx1-np2}
\end{minipage}
\end{figure}

Figure~\ref{FEx1-np1} is drawn in the case $\ell = 4$, but aims at
illustrating the possible local nodal patterns in the general case.
We look at square-like domains $\cQ_{i,j}$ which stay away from the
poles, and can be visited by the nodal set. This means that $1 \le i
\le \ell -1$, $0 \le j \le 2\ell - 1$ and that $(-1)^{i+j} = - 1$,
corresponding to the white domains on the checkerboard, see
Properties~\ref{Ex1-P1}~(i).

The local nodal pattern at the vertices is shown in the domain
labelled (A). According to Lemma~\ref{Ex1-L1}, and our assumption on
$\mu$, the nodal set $N(H^{\mu,\ell})$ consists of finitely many
disjoint simple closed regular curves. According to
Properties~\ref{Ex1-P1}~(ii), any such nodal curve can only enter a domain
$\cQ_{i,j}$ at a vertex, and exit at another one. Taking into
account Properties~\ref{Ex1-P1}~(iii), this leaves exactly three
possibilities for the nodal pattern in a domain $\cQ_{i,j}$,
illustrated in the domains labelled (B), (C) and (D). According to
the separation Lemma~\ref{Ex1-L2}, both (C) and (D) are impossible.
Notice that case (D) could also be discarded by the fact that the
nodal curves do not intersect (absence of critical zeros). Finally,
the only possible local nodal pattern in the square-like domain
$\cQ_{i,j}$ is the one shown in (B): the nodal curves ``follow'' the
meridians.

Remark. In Stern's thesis this conclusion follows from the claim
that the nodal set depends continuously on $\mu$ and that $\mu$ is
small enough.

Figure~\ref{FEx1-np2} is drawn in the case $\ell = 3$, but aims at
illustrating the possible local nodal patterns in the general case.
We look at triangle-like domains $\cQ_{i,j}$ which can be
visited by the nodal set, and one of whose vertices
is at the north or south pole.
This means that $i = 0$ or $\ell$, $0 \le j \le 2\ell - 1$ and that
$(-1)^{i+j} = - 1$, see Properties~\ref{Ex1-P1}~(i).


The same arguments as above show that there is only one possible
nodal pattern.

\subsubsection{A.~Stern's first theorem for the sphere}\label{SSS-Ex1}

We can now state the following quantitative version of A.~Stern's
first theorem, see Theorem~\ref{I-SP2odd}. Recall the notation
$\mu_c(\ell)$ in Lemma~\ref{Ex1-L1}.

\begin{proposition}\label{Ex1-P}
Assume that $0 < \mu < \mu_c(\ell)$. \vspace{-3mm}
\begin{enumerate}
    \item When $\ell$ is odd, the nodal set $N(H^{\mu,\ell})$ is a
    unique regular simple closed curve and hence, the eigenfunction
    $H^{\mu,\ell}$ has exactly two nodal domains.
    \item When $\ell$ is even, the nodal set $N(H^{\mu,\ell})$ is the union of
    $\ell$ regular disjoint simple closed curves and hence,
    the eigenfunction $H^{\mu,\ell}$ has exactly $(\ell + 1)$ nodal domains.
\end{enumerate}
\end{proposition}

\textbf{Proof}. According to the remark following
Lemma~\ref{Ex1-L1}, under the assumption on $\mu$, the nodal set of
$H^{\mu,\ell}$ is a regular $1$-dimensional submanifold. Since the
eigenfunction $H^{\mu,\ell}$ does not vanish at the north and south
poles, we can work in the exponential map at the north pole. For the
proofs below, keep in mind Section~\ref{SSS-Ex1-lnp}.

\begin{figure}
\begin{minipage}[c]{.46\linewidth}
\hspace{6mm}\includegraphics[width=6cm]{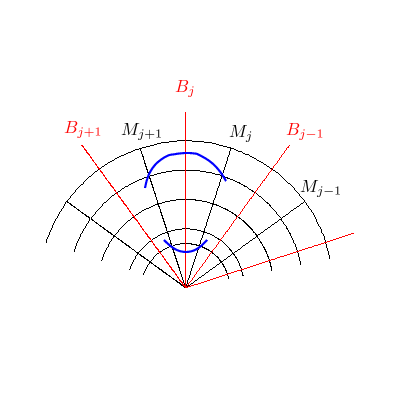}
\caption{$\ell$ even, $j$ odd}\label{Ex1-even}
\end{minipage} \hfill
\begin{minipage}[c]{.46\linewidth}
\hspace{7mm}\includegraphics[width=6cm]{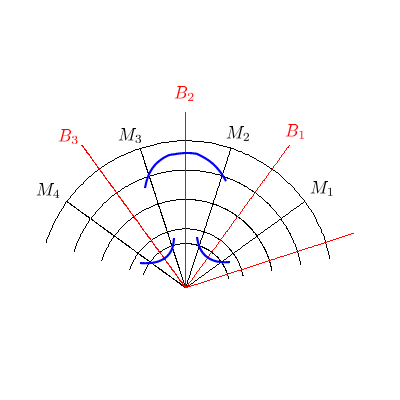}
\caption{$\ell$ odd}\label{Ex1-odd}
\end{minipage}
\end{figure}

\emph{Proof of Proposition~\ref{Ex1-P}, Assertion~(ii). The integer
$\ell$ is assumed to be even}. The proof is illustrated by
Figure~\ref{Ex1-even} which shows parts of the nodal patterns of
$Z_{\ell}$ (latitude circles) and $W_{\ell}$ (meridians) viewed in
the exponential map centered at the north pole $p_{+}$. The south
pole corresponds to the outer circle (the boundary of the maximal
domain in which the exponential map is a diffeomorphism).

Call $B'_j$ the intersection
$$B'_j = B_j \cap \{\vartheta_1(\ell) \le \vartheta \le
\pi-\vartheta_1(\ell)\}.$$ We now use Lemma~\ref{Ex1-L2}.

\noi (a)~ When $j$ is even, the function $b^{\,\mu, \ell,
j}(\vartheta)$ is positive in $[0,\vartheta_1(\ell)] \cup [\pi -
\vartheta_1(\ell), \pi]$, and hence, by Lemma~\ref{Ex1-L2},
$N(H^{\mu,\ell}) \cap B_j = \emptyset\,$. When $j$ is odd, the
function $b^{\,\mu, \ell, j}(\vartheta)$ has exactly one zero in
each of the intervals $(0,\vartheta_1(\ell))$ and $(\pi -
\vartheta_1(\ell), \pi)$. It follows that $N(H^{\mu,\ell}) \cap B_j$
consists of exactly two points, one in $\cQ_{0,j}$ and one in
$\cQ_{\ell,j}$.

\noi (b)~ Choose $j = 2k+1$, odd (see Figure~\ref{Ex1-even}). Note
that there are exactly $\ell$ such values of $j$ between $0$ and
$2\ell-1$. In $\cQ_{0,j}$ the nodal set $N(H^{\mu,\ell})$ can only
consist of a curve from the point $q_{1,j}$ to the point $q_{1,j+1}$
(see Subsection~\ref{SS-Ex1Not} for the notation, and use
Properties~\ref{Ex1-P1}), intersecting $B_j$ at exactly one point.
This curve is part of a connected component (a simple closed curve)
$\gamma_k \subset N(H^{\mu,\ell})$. We now follow the curve
$\gamma_k$, starting from $q_{1,j}$ in the direction of $q_{1,j+1}$.
According to the preceding point (a), $\gamma_k$ can meet neither
$B_{j+1}$, nor $B'_j$. Therefore, according to
Paragraph~\ref{SSS-Ex1-lnp}, the curve has to go through the points
$q_{1,j+1}, q_{2,j+1}, \ldots, q_{\ell,j+1}$ passing alternatively
inside $\cM_{j+1}$ or $\cM_j$. Since $\ell$ is even, at
$q_{\ell,j+1}$, the curve enters $\cQ_{\ell,j}$, crosses $B_j$ (at a
single point), and exits $\cQ_{\ell,j}$ at $q_{\ell,j}$. Since it
can cross neither $B_{j-1}$, nor $B'_j$, the curve $\gamma_k$ has to
go back to $q_{1,j}$, through the points $q_{\ell-1,j}, \ldots,
q_{2,j}$, alternatively inside $\cM_j$ or $\cM_{j-1}$. This means
that the simple closed curve $\gamma_k$ goes through all the points
in $\cN \cap \cM_j$, with $j = 2k+1$.

\noi (c)~ In this way, we obtain $\ell$ simple closed curves
$\gamma_1, \ldots, \gamma_{\ell}$ which are connected components of
$N(H^{\mu,\ell})$, with the curve $\gamma_k$ (where $j = 2k+1$)
contained in the sector bounded by the meridians $B_{j-1}$ and
$B_{j+1}$ and containing $B_j$. Furthermore, these $\ell$ curves
visit all the points $q_{i,j} \in \cN$. It follows from
Properties~\ref{Ex1-P1}~(iii) that there can be no other components,
and hence that
$$N(H^{\mu,\ell}) = \cup_{k=1}^{\ell} \gamma_{k}\,.$$

This finishes the proof of Proposition~\ref{Ex1-P},
Assertion~(ii).\hfill $\square$ \medskip

\emph{Proof of Proposition~\ref{Ex1-P}, Assertion~(i). The integer
$\ell$ is now assumed to be odd}. The proof is illustrated by
Figure~\ref{Ex1-odd} which shows parts of the nodal patterns of
$Z_{\ell}$ (latitude circles) and $W_{\ell}$ (meridians) viewed in
the exponential map centered at the north pole $p_{+}\,$. The south
pole corresponds to the outer circle (the boundary of the domain in
which the exponential map is a diffeomorphism).

As in the previous proof, call $B'_j$ the intersection
$$B'_j = B_j \cap \{\vartheta_1(\ell) \le \vartheta \le
\pi-\vartheta_1(\ell)\}.$$

\noi (a)~ When $j$ is even, the function $b^{\,\mu, \ell,
j}(\vartheta)$ is positive in $[0,\vartheta_1(\ell)]$, and admits exactly
one zero in the interval $(\pi - \vartheta_1(\ell), \pi)$. Hence
$N(H^{\mu,\ell}) \cap B_j$ contains exactly one point located in
$\cQ_{\ell,j}$. When $j$ is odd, the function $b^{\,\mu, \ell,
j}(\vartheta)$ has exactly one zero in the interval
$\left( 0,\vartheta_1(\ell) \right)$ and is negative in $[\pi - \vartheta_1(\ell), \pi]$. By
Lemma~\ref{Ex1-L2}, it follows that $N(H^{\mu,\ell}) \cap B_j$
consists of exactly one point located in $\cQ_{0,j}$.

\noi (b)~ Choose $j = 1$. In $\cQ_{0,1}$ the nodal set
$N(H^{\mu,\ell})$ can only consist of a curve going from the point
$q_{1,1}$ to the point $q_{1,2}$ (see Subsection~\ref{SS-Ex1Not} for
the notation and use Properties~\ref{Ex1-P1}), intersecting $B_1$
at exactly one point. This curve is part of a connected component (a
simple closed curve) $\gamma \subset N(H^{\mu,\ell})$. We now follow
the curve $\gamma$, starting from $q_{1,1}$ in the direction of
$q_{1,2}$. According to the preceding point (a), $\gamma$ can meet
neither $B'_{2}$, nor $B_1$, so that it has to go through the points
$q_{1,2}, q_{2,2}, \ldots, q_{\ell,2}$ passing alternatively inside
$\cM_{2}$ or $\cM_1$. Because $\ell$ is odd, at $q_{\ell,2}$, the
curve exits $\cQ_{\ell-1,1}$, enters $\cQ_{\ell,2}$, crosses $B_2$
(at a single point), and exits $\cQ_{\ell,2}$ at $q_{\ell,3}$ into
$\cM_3$. Since it can cross neither $B'_{3}$, nor $B'_2$, the curve
$\gamma$ has to go to $q_{1,3}$, through the points $q_{\ell-1,3},
\ldots q_{2,3}$, alternatively inside $\cM_2$ or $\cM_{3}$. The
curve $\gamma$ therefore goes from $q_{1,1}$ to $q_{1,3}$ where we
can start again with the same argument as before. Iterating $\ell$
times the argument, the curve $\gamma$ gets back to its initial
point $q_{1,1}$.

\noi (c)~ In this way, we obtain a simple closed curve $\gamma$ in
$N(H^{\mu,\ell})$, which crosses all the meridians $B_j$ once, and
which visits all the points $q_{i,j} \in \cN$. It follows from
Properties~\ref{Ex1-P1}~(iii) that there can be no other component,
and hence that
$$N(H^{\mu,\ell}) = \gamma\,.$$

This finishes the proof of Proposition~\ref{Ex1-P},
Assertion~(i).\hfill $\square$ \medskip

\begin{figure}[!ht]
\begin{center}
\includegraphics[width=12cm]{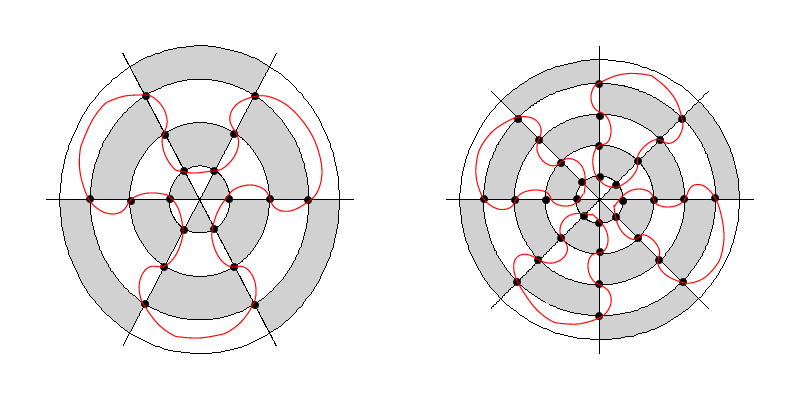}
\caption{Cases $\ell = 3$ and $\ell = 4$} \label{FEx1-nodalset}
\end{center}
\end{figure}

It is easy to follow the above proofs on Figure~\ref{FEx1-nodalset}
which shows the nodal set of $H^{\mu,\ell}$, in the exponential map
at $p_{+}$, for $\mu > 0$ small enough and for $\ell = 3$ (left) and
$\ell = 4$ (right).


\section{Stern's second theorem: even case}\label{S-Ex2}

The purpose of this section is to prove Theorem~\ref{I-SP2even}. As
a matter of fact, we shall give a more quantitative result,
Proposition~\ref{Ex2-P}, which implies the theorem. As in
Section~\ref{S-Ex1}, we follow the ideas of A.~Stern sketched in the
introduction.

Fix an integer $$\ell = 2r \ge 2\,,$$
 as well as an angle $\alpha$ defined by
\begin{equation}\label{defalpha}
\alpha =
\frac{\epsilon \pi}{2r}\,, \mbox{ with } 0 < \epsilon < \frac{1}{2}\,.
\end{equation}

\subsection{Notation}\label{SS-Ex2Not}

As in the first example, we consider the spherical harmonic of
degree $\ell = 2r$
\begin{equation}\label{Ex2-a1}
W(x,y,z) = \Im (x+iy)^{2r}\,,
\end{equation}
whose expression in spherical coordinates $(\vartheta,\varphi) \in\,
(0,\pi)\times\R_{2\pi}$ is given by
\begin{equation}\label{Ex2-a2}
w(\vartheta,\varphi) = \sin^{2r}(\vartheta)\, \sin(2r\varphi)\,.
\end{equation}

The perturbation of $W$ is chosen to be the spherical harmonic
$V_{\alpha}$, of degree $2r$, whose expression in spherical
coordinates is given by
\begin{equation}\label{Ex2-a3}
v_{\alpha}(\vartheta,\varphi) = P^1_{2r}(\cos \vartheta) \,
\sin(\varphi-\alpha)\,.
\end{equation}
According to Properties~\ref{P-P1}~(i), we have $P^{1}_{2r}(t) = -
(1-t^2)^{1/2} \, \frac{d}{dt}P_{2r}(t)$, so that
\begin{equation}\label{Ex2-a4}
v_{\alpha}(\vartheta,\varphi) = - \sin \vartheta \, P'_{2r}(\cos
\vartheta) \, \sin(\varphi-\alpha)\,.
\end{equation}
According to \eqref{SH-1}, the corresponding harmonic homogeneous
polynomial of degree $2r$ in $\R^3$ is given by the formula
\begin{equation}\label{Ex2-a5}
V_{\alpha}(x,y,z) = \left( \sin\alpha \,x - \cos\alpha \,y\right)
\sum_{j=0}^{r-1}a_j z^{2r-2j-1}(x^2+y^2+z^2)^{j}\,,
\end{equation}
where the $a_j$'s are the coefficients of the polynomial $P'_{2r}$,
$$
P'_{2r}(t) = \sum_{j=0}^{r-1} a_j t^{2r-2j-1}\,.
$$\smallskip

The nodal set of the spherical harmonic $W$ consists of the $2\ell =
4r$ meridians $M_j, 0 \le j \le 4r-1\,$, defined as in
Subsection~\ref{SS-Ex1Not},
$$
N(W) = \bigcup_{j=0}^{4r-1} M_j \,,
$$
with the corresponding open sectors $\cM_j$ on the sphere.

These meridians meet at the north and south poles which are the only
critical zeros of $W$, $W(p_{\pm})=0$, and $d_{p_{\pm}}W = 0\,$ (the
differential of the function $W$ at the poles).
\medskip

The nodal set of the spherical harmonic $V_{\alpha}$ consists of
$(2r-1)$ latitude circles $L'_i$ ($1 \le i \le 2r-1$), and two
meridians $M'_0$ and $M'_1$,
\begin{equation}\label{Ex2-b1}
N(V_{\alpha}) = \bigcup_{j=0}^{2r-1} L'_i \, \bigcup M'_0 \bigcup
M'_1 \,.
\end{equation}

The latitude circles,
\begin{equation}\label{Ex2-b2}
L'_i = \{(\vartheta,\varphi) ~|~ \vartheta = \vartheta'_i(2r)\}\,,
~1\le i \le 2r-1\,,
\end{equation}
are associated with the $(2r-1)$ zeros, $~0 < \vartheta'_1(2r) <
\cdots < \vartheta'_{2r-1}(2r) < \pi~$, of the function
$P'_{2r}(\cos\vartheta)$, see Properties~\ref{P-P1}~(iii), and we
let $\vartheta'_0(2r)=0$ and $\vartheta'_{2r}(2r)=\pi$. They
determine sectors
\begin{equation}\label{Ex2-b4}
\cL'_i = \{(\vartheta,\varphi) ~|~ \vartheta'_i(2r) < \vartheta <
\vartheta'_{i+1}(2r) \}\,, 0\le i \le 2r-1 \,.
\end{equation}

The meridians $M'_{k}$ are given by
\begin{equation}\label{Ex2-b6}
M'_0 = \{(\vartheta,\varphi) ~|~ \varphi = \alpha\} \text{~and~}
M'_1 = \{(\vartheta,\varphi) ~|~ \varphi = \alpha + \pi\} \,.
\end{equation}
They determine sectors
\begin{equation}\label{Ex2-b8}
\cM'_0 = \{(\vartheta,\varphi) ~|~ \alpha < \varphi < \alpha + \pi\}
\text{~and~} \cM'_1 =\{(\vartheta,\varphi) ~|~ \alpha + \pi <
\varphi < 2\pi + \alpha \}\,.
\end{equation}

\begin{figure}
\begin{minipage}[c]{.46\linewidth}
\includegraphics[width=6cm]{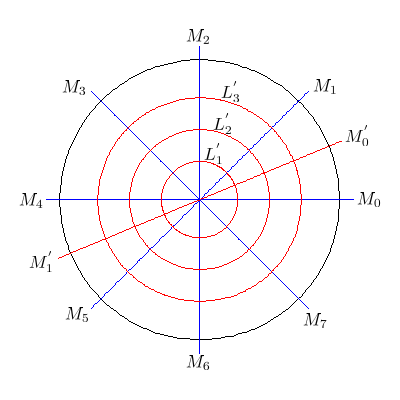}
\caption{$N(W)$ and $N(V_{\alpha})$, $\ell = 4$}
\label{FEx2-nodal1}
\end{minipage} \hfill
\begin{minipage}[c]{.46\linewidth}
\includegraphics[width=6cm]{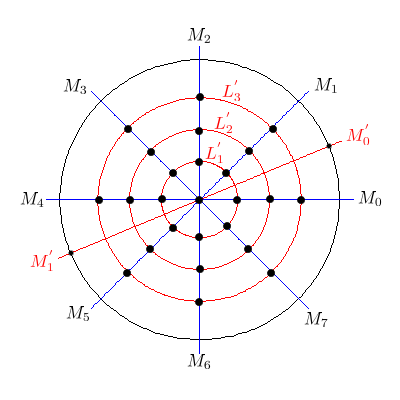}
\caption{$\cN$, $\ell = 4$} \label{FEx2-N}
\end{minipage}
\end{figure}

Figure~\ref{FEx2-nodal1} shows the nodal sets $N(W)$ and
$N(V_{\alpha})$, in the case $\ell = 2r = 4\,$. They are viewed in the
exponential map $\exp_{p_{+}}$\, \ie, in the disk $D(0,\pi)$, whose
boundary corresponds to the cut-locus of $p_{+}$ \ie, $p_{-}$\,.

As in the first example, the set $\cN = N(W) \cap N(V_{\alpha})$ of
common zeros to the spherical harmonics $W$ and $V_{\alpha}$, plays
a special role. We have

\begin{equation}\label{Ex2-b10}
\cN = \{p_{+},p_{-}\} \bigcup \{q_{i,j} ~|~ 1 \le i \le 2r-1, ~0 \le
j \le 4r-1\}\,,
\end{equation}
where $q_{i,j}$ is the intersection point of the latitude circle
$L'_i$ with the meridian $M_j$.\\ Figure~\ref{FEx2-N} shows the set
$\cN$ in the exponential map. The points in $\cN$ appear as the big
dots: the intersection points of the latitude circles $L'_i$, with
the meridians $M_j$, and the poles. Note that the south pole is
represented by two small dots, one on the meridian $M'_0$, one on
the meridian $M'_1$. Note that there are no other dots on these
meridians, see Properties~\ref{Ex2-P1}~(iii).


We also introduce the connected components of the set $\{W \,
V_{\alpha} \neq 0\}$,
\begin{equation}\label{Ex2-b12}
\cQ_{i,j,k} = \cL'_i \cap \cM_j \cap \cM'_k \,, ~0\le i\le 2r-1,
~0\le j \le 4r-1, ~k = 0, 1 \,.
\end{equation}
Note that
\begin{equation}\label{Ex2-b14}
\mathrm{sgn}(W\,V_{\alpha}) = (-1)^{i+j+k+1} \text{~on~} \cQ_{i,j,k}
\,.
\end{equation}

\subsection{The family $H^{\mu}$}\label{SS-Ex2Fam}

\subsubsection{Definition}
Following the ideas of Stern \cite{St}, we consider the
one-parameter family of spherical harmonics of degree $\ell = 2r$,
\begin{equation}\label{Ex2-c2}
H^{\mu}(x,y,z) = W(x,y,z) - \mu V_{\alpha}(x,y,z)\,,
\end{equation}
whose expression in spherical coordinates is given by
\begin{equation}\label{Ex2-c4}
\begin{split}
h^{\mu}(\vartheta,\varphi) & = w(\vartheta,\varphi) - \mu \,
v_{\alpha}(\vartheta,\varphi)\,,\\
& = \sin^{2r}(\vartheta)\,\sin(2r\varphi) + \mu \, \sin\vartheta \,
P'_{2r}(\cos\vartheta) \, \sin(\varphi -\alpha)\,.
\end{split}
\end{equation}

Note that
\begin{equation}\label{Ex2-sgnmu}
h^{-\mu}(\vartheta,\varphi) = h^{\mu}(\vartheta,\varphi+\pi)\,.
\end{equation}
It follows that it suffices to consider the case $\mu > 0$. We shall
therefore assume that $\mu > 0$ for the remainder of
Section~\ref{S-Ex2}.

\subsubsection{Critical zeros}\label{SSS-Ex2-cz}

We now investigate the critical zeros of $H^{\mu}$. The spherical
harmonic $W$ vanishes at order at least $3$ at the poles, while the
nodal set of $V_{\alpha}$ is a piece of great circle at each pole.
It follows that the north and south poles are not critical points of
$N(H^{\mu})$, see also Properties~\ref{Ex2-P1}~(i). We can therefore
look for critical zeros in the spherical coordinates \ie, look for
critical zeros of $h^{\mu}$ in $(0,\pi)\times \R_{2\pi}$.

The point $(\vartheta,\varphi) \in\, (0,\pi)\times \R_{2\pi}$ is a
critical zero of $h^{\mu}$ if and only if,
\begin{equation}\label{Ex2-d2}
\begin{split}
& h^{\mu}(\vartheta,\varphi) = 0 \,,\\
& \partial_{\vartheta} h^{\mu}(\vartheta,\varphi) = 0 \,,\\
& \partial_{\varphi} h^{\mu}(\vartheta,\varphi) = 0 \,.\\
\end{split}
\end{equation}

Using the second order differential equation \eqref{P-2} satisfied
by the Legendre polynomial $P_{2r}$, we find that
\begin{equation*}
\begin{split}
\partial_{\vartheta} h^{\mu}(\vartheta,\varphi) = & ~2r \, \cos\vartheta
\, \sin^{2r-1}(\vartheta)\, \sin(2r\varphi)\\
& ~ + \mu \sin(\varphi - \alpha) \big[ 2r(2r+1)
P_{2r}(\cos\vartheta) - \cos\vartheta \, P'_{2r}(\cos\vartheta)\big]
\,.
\end{split}
\end{equation*}

It follows that the point $(\vartheta,\varphi) \in\, (0,\pi)\times
\R_{2\pi}$ is a critical zero of $h^{\mu}$ if and only if,
\begin{equation}\label{Ex2-d4}
\begin{split}
& \sin^{2r-1}(\vartheta) \, \sin(2r\varphi) + \mu \,
P'_{2r}(\cos\vartheta)\, \sin(\varphi -\alpha) = 0 \,,\\[4pt]
& 2r \, \cos\vartheta \, \sin^{2r-1}(\vartheta)\, \sin(2r\varphi)\\
&  ~~+ \mu \sin(\varphi - \alpha) \big[ 2r(2r+1)
P_{2r}(\cos\vartheta) - \cos\vartheta P'_{2r}(\cos\vartheta)\big] = 0 \,,\\[4pt]
& 2r \, \sin^{2r-1}(\vartheta) \, \cos(2r\varphi) + \mu \,
P'_{2r}(\cos \vartheta)\,
\cos(\varphi -\alpha) = 0 \,.\\
\end{split}
\end{equation}

The pair of the first and third equations in \eqref{Ex2-d4} is equivalent
to the pair of the first and third equations  in \eqref{Ex2-d6} below. Plugging
the first equation in \eqref{Ex2-d4} into the second one, and using the
fact that $\mu > 0$, we get the
second equation in \eqref{Ex2-d6}. It follows that the point
$(\vartheta,\varphi) \in\, (0,\pi)\times \R_{2\pi}$ is a critical
zero of $h^{\mu}$ if and only if,
\begin{equation}\label{Ex2-d6}
\begin{split}
& \mu \, P'_{2r}(\cos\vartheta) + \sin^{2r-1}(\vartheta)\,\big[ 2r
\cos(2r\varphi)\,\cos(\varphi -\alpha) +
\sin(2r\varphi)\,\sin(\varphi -\alpha) \big] = 0 \,,\\
& \sin(\varphi -\alpha) \, \big[ 2r P_{2r}(\cos\vartheta) -
\cos\vartheta P'_{2r}(\cos\vartheta)\big]  = 0 \,,\\
& 2r \cos(2r\varphi)\,\sin(\varphi -\alpha) -
\sin(2r\varphi)\,\cos(\varphi -\alpha)  = 0 \,.
\end{split}
\end{equation}

\begin{property}\label{Ex2-P0}
Assume that $\mu > 0$. Then, the product $\sin(2r\varphi)\,
\sin(\varphi - \alpha)$ does not vanish at the critical zeros of
$H^{\mu}$.
\end{property}

This property follows from the third equation in
\eqref{Ex2-d6}.\medskip

Finally, it follows that the point $(\vartheta,\varphi) \in\,
(0,\pi)\times \R_{2\pi}$ is a critical zero of $h^{\mu}$ if and only
if,
\begin{equation}\label{Ex2-d8}
\begin{split}
\mu \, P'_{2r}(\cos\vartheta) + \sin^{2r-1}(\vartheta)\,\big[ 2r
\cos(2r\varphi)\,\cos(\varphi -\alpha) +
\sin(2r\varphi)\,\sin(\varphi -\alpha) \big] & = 0 \,,\\
2r P_{2r}(\cos\vartheta) - \cos\vartheta P'_{2r}(\cos\vartheta) & = 0 \,,\\
2r \cos(2r\varphi)\,\sin(\varphi -\alpha) -
\sin(2r\varphi)\, \cos(\varphi -\alpha) & = 0 \,.\\
\end{split}
\end{equation}

We first analyze the second equation in \eqref{Ex2-d8}. Define $Q(t)
:= 2r P_{2r}(t) - t P'_{2r}(t)\,$. This is an even polynomial of
degree less than or equal to $(2r-2)$. For parity reasons the roots
of the polynomials $Q, P_{2r}$ and $P'_{2r}$ are symmetric with
respect to $0\,$, and it suffices to look at $t\ge0\,$. According to
Properties~\ref{P-P1}~(iii), the non-negative roots $t_i$ of $P_{2r}\,$,
and $t'_i$ of $P'_{2r}$ satisfy
\begin{equation*}
0 =t'_r < t_r < t'_{r-1} < t_{r-1} < \cdots < t_2 < t'_1 < t_1 < 1\,.
\end{equation*}
The following equalities are easy to check,
\begin{equation*}
\sgn{P_{2r}(t'_i)} = (-1)^{i} \text{~and~} \sgn{P'_{2r}(t_i)} =
(-1)^{i-1}\,,
\end{equation*}
\begin{equation*}
\sgn{Q(t_i)} = (-1)^{i} \text{~and~} \sgn{Q(t'_{i-1})} =
(-1)^{i-1}\,.
\end{equation*}
It follows that $Q$ vanishes at least once in each interval
$(t_{i+1},t'_{i})$ for $1\le i \le r-1\,$. Since $Q$ has at most
$(r-1)$ non-negative zeros, we can conclude that $Q$ has exactly
$(r-1)$ zeros in $(0,1)$, and more precisely one zero, which we
denote by $\cos\omega_{i}$, in each interval $(t_{i+1},t'_{i})$, so
that $\omega_{i} \in \, (\vartheta'_{i}(2r),\vartheta_{i+1}(2r))$,
and
\begin{equation*}
0 < \vartheta_1(2r) < \vartheta'_1(2r) < \omega_1 < \vartheta_2(2r)
< \cdots < \vartheta'_{r-1}(2r) < \omega_{r-1} < \vartheta_r(2r) <
\vartheta'_{r}(2r) = \frac{\pi}{2}\,.
\end{equation*}
Note that the inequalities are strict \ie, that
$P_{2r}(\cos\omega_i) \neq 0$ and $P'_{2r}(\cos\omega_i) \neq 0\,$,
and that the zeros $\omega_i$ depend on $\ell = 2r\,$. \medskip

We now analyze the third equation in \eqref{Ex2-d8}. Define the
function,
\begin{equation}\label{Ex2-d9}
f(\varphi) = 2r \, \cos(2r\varphi)\, \sin(\varphi -\alpha) -
\sin(2r\varphi) \, \cos(\varphi - \alpha)\,.
\end{equation}
The function $f$ satisfies $f(\pi + \varphi) + f(\varphi) = 0\,$, and
$f'(\varphi) = - (4r^2-1) \sin(2r\varphi) \, \sin(\varphi -\alpha)\,$.
 An easy analysis in $[0,\pi]$ (using the choice of $\alpha$)
shows that $f$ does not vanish in $[0,\frac{\pi}{2r}]$, and has
exactly one zero in each interval
$[j\frac{\pi}{2r},(j+1)\frac{\pi}{2r}]$ for $1 \le j \le 2r-1\,$. It
follows that $f$ has exactly $(4r-2)$ zeros in $[0,2\pi]$, $0 <
\varphi_1 < \varphi_2 < \cdots < \varphi_{4r-2} < 2\pi$, and that
$\varphi_j \in \left(j\frac{\pi}{2r},(j+1)\frac{\pi}{2r}\right)$.

The only possible critical zeros of $H^{\mu}$ are given in spherical
coordinates by the points  $(\omega_i,\varphi_j)$ and $(\pi -
\omega_i,\varphi_j)\,$, for $1 \le i \le r-1$ and $1 \le j \le
4r-2\,$. These points can only occur as critical zeros for finitely
many values of $\mu$ given by the first equation in \eqref{Ex2-d8}.
Since we work with $\mu > 0$, these critical values of $\mu$ are
given (see also the first line in \eqref{Ex2-d4}) by
\begin{equation}\label{Ex2-muv}
\mu_{i,j}(\alpha) = \Big| \frac{\sin^{2r-1}(\omega_i) \,
\sin(2r\varphi_j)}{P'_{2r}(\cos\omega_i)\,\sin(\varphi_j-\alpha)}
\Big| \,,
\end{equation}
for $1\le i \le r-1, ~~ 1 \le j \le 2r-1\,$, and can be numerically
computed.\\
We summarize the preceding analysis  in the following lemma. \medskip

\begin{lemma}\label{Ex2-L1} For $\mu > 0\,$, the spherical
harmonic $H^{\mu}$ has no critical zero except for finitely many
values of $\mu$ which are given by \eqref{Ex2-muv}. For each
value $\mu_{i,j}(\alpha)\,$, the spherical harmonic $H^{\mu}$ has
finitely many critical zeros. Define the number $\mu_{c}(\alpha,2r)$
to be
\begin{equation}\label{Ex2-mua}
\mu_{c}(\alpha,2r) = \inf \mu_{i,j}(\alpha) \,,
\end{equation}
where the infimum is taken over $1 \le i \le r-1\,$, and $1 \le j \le
2r-1\,$. Then, for $0 < \mu < \mu_{c}(\alpha,2r)\,$, the function
$H^{\mu}$ has no critical zero, so that its nodal set $N(H^{\mu})$
consists of finitely many disjoint simple closed curves.
\end{lemma}

\textbf{Remark}. One can bound  $\mu_{c}(\alpha,2r)$  from
below using the inequalities satisfied by the $\omega_i$, and
Properties~\ref{P-P1}~(vii).

\subsubsection{A separation lemma for $N(H^{\mu})$}\label{SSS-Ex2-sl}

 For $1 \le j \le 4r-2\,$, call $C_j$ the meridian,
\begin{equation}\label{Ex2-MC}
C_j = \{(\vartheta,\varphi) ~|~ \varphi = \varphi_j\}\,,
\end{equation}
where $\varphi_j$ are the zeros of the function $f$ defined in
\eqref{Ex2-d9}.

We now look at the restriction of the spherical
harmonic $H^{\mu}$ to the meridian $C_j$. Recall from
Properties~\ref{P-P1}~(iii), that $\cos\vartheta'_1(2r)$ is the
largest zero of the function $P'_{2r}(\cos\vartheta)$ in $[0,\pi]$.

\begin{lemma}\label{Ex2-L2}
Define the functions,
\begin{equation}\label{Ex2-c8}
\begin{split} b^{\,\mu,j}(\vartheta) & =
h^{\mu}(\vartheta,\varphi_j) \\
& = \sin(2r\varphi_j)\,\sin^{2r}(\vartheta) + \mu\, \sin\vartheta \,
P'_{2r}(\cos\vartheta) \,\sin(\varphi_j - \alpha)\,.
\end{split}
\end{equation}
Assume that $2r+1 \le j \le 4r-1$.\vspace{-3mm}
\begin{enumerate}
    \item For $0 < \mu < \mu_{c}(\alpha,2r)$, the functions $b^{\,\mu,j}(\vartheta)$ do not
    vanish in the interval $[\vartheta'_1(2r),\pi-\vartheta'_1(2r)]$.
    \item When $j$ is odd, the function $b^{\,\mu,j}(\vartheta)$ vanishes exactly
    once in the interval $\left(0,\vartheta'_1(2r)\right)$, and does not
    vanish in the interval $[\pi-\vartheta'_1(2r),\pi]$.
    \item When $j$ is even, the function $b^{\,\mu,j}(\vartheta)$ vanishes exactly
    once in the interval  $(\pi-\vartheta'_1(2r),\pi)$, and does not
    vanish in the interval $[0,\vartheta'_1(2r)]$.
\end{enumerate}\vspace{-3mm}
The above assertions determine the possible intersections of the
nodal set $N(H^{\mu})$ with the meridian $C_j$.
\end{lemma}

\textbf{Proof}. { Notice that the assumptions on $j$ and
$\alpha$ imply that $\sin(\varphi_j-\alpha) < 0$, and that $(-1)^j
\sin(2r\varphi_j) > 0$. The function $b^{\,\mu,j}$ vanishes at the
points such that
$$
\frac{\sin\vartheta \, P'_{2r}(\cos\vartheta)}{\sin^{2r}(\vartheta)}
= - \frac{\sin(2r\varphi_j)}{\mu\, \sin(\varphi_j-\alpha)}\,.
$$
Call $\beta_2(\vartheta)$ the function in the left-hand side. Using
the differential equation satisfied by $P_{2r}$, one finds that
$$
\beta'_2(\vartheta) = (2r+1)\, \frac{2r\,
P_{2r}(\cos\vartheta)-\cos\vartheta \,
P'_{2r}(\cos\vartheta)}{\sin^{2r}(\vartheta)}\,.
$$

The local extrema of $\beta_2$ in the interval
$[\vartheta'_1(2r),\pi-\vartheta'_1(2r)]$ are achieved at the zeros
$\omega_i \in \left(\vartheta'_i(2r),\vartheta_{i+1}(2r) \right)$ of
the second equation in \eqref{Ex2-d8}, for $1 \le i \le 2r-2$. We
have at these points,
$$
h^{\,\mu}(\omega_i,\varphi_j) = \sin(2r\varphi_j) \,
\sin^{2r}(\omega_i) \,\Big[ 1 + \mu \, \frac{P'_{2r}(\cos\omega_i)\,
\sin(\varphi_j-\alpha)}{\sin(2r\varphi_j)\,\sin^{2r-1}(\omega_i)}\Big].
$$
The first assertion follows from \eqref{Ex2-muv}.}

We now  determine what happens in the intervals
$\left(0,\vartheta'_1(2r)\right)$ and $(\pi-\vartheta'_1(2r), \pi)$.

Write $b^{\,\mu,j}(\vartheta) = (-1)^j\, \sin\vartheta \,
f_j(\vartheta)$. The derivative $f'_j(\vartheta)$ is given by
$$
f'_j(\vartheta) = (2r-1) \, (-1)^j \, \sin(2r\varphi_j) \,
\cos\vartheta \, \sin^{2r-2}(\vartheta) - (-1)^j \, \mu \,
\sin\vartheta \, P^{''}_{2r}(\cos\vartheta)
\sin(\varphi_j-\alpha)\,.
$$

Recall that $P_{2r}$ and $P^{''}_{2r}$ are even functions and that
$P'_{2r}$ is odd. The largest zeros of theses functions in $[-1,1]$
satisfy, with an obvious notation,
$$
 t_2 < t^{''}_1 < t'_1 < t_1 < 1 \,.
$$
Looking at the signs of these functions in the various intervals
between $ t_2$ and $1$, and using the parity to determine what
happens near $-1\,$, we can make the following observations.

\noib Case $j$ even. For $\vartheta \in \left( 0,\vartheta'_1(2r) \right)$,
$f_j(\vartheta) > 0$. On the other-hand, $f_j\left( \pi-\vartheta'_1(2r) \right) > 0$
and $f_j(\pi) < 0\,$, while $f'_j(\vartheta) < 0$ in
$\left( \pi-\vartheta'_1(2r),\pi \right)\,$.

\noib Case $j$ odd. For $\vartheta \in \left( \pi-\vartheta'_1(2r),\pi \right)$,
$f_j(\vartheta) > 0\,$. On the other-hand, $f_j(0) < 0$ and
$f_j\left( \vartheta'_1(2r) \right) > 0\,$, while $f'_j(\vartheta) >$ in $\left( 0,
\vartheta'_1(2r) \right)\,$.

The second and third assertion follows. \hfill $\square$\medskip

\subsubsection{General properties of $N(H^{\mu})$}\label{SSS-Ex2-gp}

\begin{properties}\label{Ex2-P1}
For $\mu > 0$, the nodal sets $N(H^{\mu})$ share the following
properties.\vspace{-3mm}
\begin{enumerate}
    \item The north and south poles are zeros of order $1$ of $H^{\mu}$,
    $H^{\mu}(p_{\pm}) = 0$, and $d_{p_{\pm}}H^{\mu} \neq 0$. In particular,
    near each pole, the nodal set $N(H^{\mu})$ consists of a single
    arc, tangent to the great circle $M'_0 \cup M'_1$.
    \item The nodal set $N(H^{\mu})$ satisfies,
    \begin{equation}\label{Ex2-c6}
    \cN \subset N(H^{\mu}) \subset \cN \, \bigcup\, \{W \,
    V_{\alpha} > 0\} \,.
    \end{equation}
    \item Since $\alpha = \frac{\epsilon \pi}{2r}$, with $0 < \epsilon
    < \frac{1}{2}$, the nodal set $N(H^{\mu})$ meets the great circle $M'_0 \cup M'_1$
    at the poles tangentially, and nowhere else.
    \item The connected components of $N(H^{\mu})$ are contained in
    either the closed hemisphere $\overline{\cM'_0}$ or in
    $\overline{\cM'_1}$.
        \item The points in $\cN$ (common zeros to the spherical
    harmonics $W$ and $V_{\alpha}$) are not critical zeros of $H^{\mu}$.
    At the points $q_{i,j} \in \cN, ~1 \le i \le 2r-1, ~0 \le j \le
    4j-1\,$, the nodal set $N(H^{\mu})$ consists of a single arc
    which is transversal to the latitude circles $L'_i$ and to the
    meridians $M_j$.
    \item For $0 < \mu < \mu_{c}(\alpha,2r)$ (defined in \eqref{Ex2-mua}), no closed
    component of the nodal
    set $N(H^{\mu})$ can be entirely contained in some domain $\cQ_{i,j,k}$.
\end{enumerate}
\end{properties}

\textbf{Proof}. Assertion~(i) follows from the fact that near the
poles, $N(V_{\alpha})$ is a piece of great circle, while $W$
vanishes at order at least $3$. Assertion~(ii) is clear (this is the
checkerboard property introduced by A.~Stern as we recalled in the
introduction). Assertion~(iii) is clear because the great circle
$M'_0 \cup M'_1$ only meets the nodal set $N(W)$ at the poles.
Assertion~(iv) follows from Assertion~(ii), the choice of $\alpha$
and the parity of $\ell = 2r$. We can indeed look at a neighborhood of the north
pole (the pattern near the south pole is the image of the pattern at
the north pole under the antipodal map). The nodal curve at $p_{+}$
must visit the domains $\cQ_{0,0,1}$ and $\cQ_{0,2r,1}$, both in
$\cM'_1$, and cannot visit the domains $\cQ_{0,0,0}$ and
$\cQ_{0,2r,0}$. On the other-hand, as we already pointed out, the
nodal set cannot meet the great circle $M'_0 \cup M'_1$.
Assertion~(v) follows by checking that the partial derivatives
$\partial_{\vartheta} h^{\mu}$ and $\partial_{\varphi} h^{\mu}$ do
not vanish at the points $\left( \vartheta'_i(2r),\frac{j\pi}{2r}
\right)$.
Assertion~(vi) follows by using the same energy argument as in
Properties~\ref{Ex1-P1}. \hfill $\square$\medskip

Figure~\ref{FEx2-nodal2} illustrates the proof of
Properties~\ref{Ex2-P1}. The checkerboard appears in white/grey
(allowed / forbid\-den domains).


\subsubsection{Local nodal patterns for $H^{\mu}$}\label{SSS-Ex2-lnp}

The arguments to determine the local nodal patterns for $H^{\mu}$
are the same as in Paragraph~\ref{SSS-Ex1-lnp},  with an extra case.
Namely, at each pole, the nodal set $N(H^{\mu})$ is a single arc
tangent to the great circle $M'_0 \cup M'_1$, going through two
triangle-like domains, one of whose vertices is the pole. One of the
two remaining vertices does not belong to $\cN$, the other does
belong to $\cN$ so that the local nodal pattern is well determined.
See Figures~\ref{FEx1-np1}, \ref{FEx1-np2} and \ref{FEx2-nodal2}.

\subsubsection{A.~Stern's second theorem}\label{SSS-Ex2n}

We can now state the following improved version of Stern's second
theorem, Theorem~\ref{I-SP2even}. Recall the definition of
$\mu_c(\alpha,2r)$ given  in \eqref{Ex2-mua}.

\begin{proposition}\label{Ex2-P}
For $\alpha$ satisfying \eqref{defalpha} and  $0 < \mu <
\mu_c (\alpha,2r)\,$, \vspace{-3mm}
\begin{enumerate}
    \item the spherical harmonic $H^{\mu}$, of degree $2r$,  introduced in \eqref{Ex2-c2},
    has no critical zero,
    \item the nodal set $N(H^{\mu})$ of $H^{\mu}$ has exactly two
    connected components \ie, consists of exactly two simple closed
    curves which do not intersect.
\end{enumerate} \vspace{-3mm}
In particular, for $0 < \mu < \mu_c(\alpha,2r)$, the spherical harmonic
$H^{\mu}$ has exactly three nodal domains.
\end{proposition}

\textbf{Proof of Proposition~\ref{Ex2-P}}

Note that $H^{\mu}$ is even, so that it is invariant under the
antipodal map, and so is its nodal set $N(H^{\mu})$. We have already
seen, Properties~\ref{Ex2-P1}, that a connected component of
$N(H^{\mu})$ is contained in either $\overline{\cM'_0}$ or
$\overline{\cM'_1}$. Furthermore, there is one connected component,
call it $\gamma$, which is contained in $\overline{\cM'_1}$, and
which is tangent to the great circle $M'_0 \cup M'_1$ at the north
pole $p_{+}$. Similarly, there is another connected component which
is contained in $\overline{\cM'_0}$, and which is tangent to $M'_0
\cup M'_1$ at the south pole $p_{-}$. The second can be deduced from
$\gamma$ by applying the antipodal map.

It follows that it suffices to look at the part of the nodal set
$N(H^{\mu})$ which is contained in $\overline{\cM'_1}$.  For
this reason, we only have to consider the meridians $C_j$ for $2r+1
\le j \le 4r-1$. The connected component $\gamma$ is a simple
closed curve. Start from the north pole, tangentially to $M'_0$,
inside the domain $\cQ_{0,0,1}\,$. The only possibility for $\gamma$
is to exit $\cQ_{0,0,1}$ through the point $q_{1,0,1}$. Using the
separation lemma, Lemma~\ref{Ex2-L1}, and the analysis of local
nodal patterns, we see that $\gamma$ has to wind around $M_0$,
inside the white domains, until it reaches the last point
$q_{2r-1,0,1}\,$, at which it has to enter the white domain
$\cQ_{2r-1,4r-1,1}\,$, cross the meridian $B_{4r-1}$, exit through
the point $q_{2r-1,4r-1,1}$ and wind along $M_{4r-1}$ until it
reaches the domain $\cQ_{0,4r-2,1}\,$, \etc. The situation is
similar to the one we encountered in the proof of
Proposition~\ref{Ex1-P}~(i). Indeed, the important point in this
proof was that the number $\ell$ of latitude circle $L_i$ was odd.
In the present case we have $\ell = 2r$, but the number of latitude
circles $L'_i$ is $2r-1$, an odd integer. The proof of
Section~\ref{S-Ex1} applies mutatis mutandis, and the conclusion is
that $\gamma$ goes back to the north pole after going up and down
$r$ times, visiting all the points in $\cN \cap \cM'_1$. Using
Properties~\ref{Ex2-P1}, it follows that $N(H^{\mu}) \cap
\overline{\cM'_1}$ has exactly one connected component $\gamma$.
Using the antipodal map, this means that $N(H^{\mu})$ has exactly
two connected components. \hfill $\square$

\begin{figure}
\begin{minipage}[c]{.46\linewidth}
\includegraphics[width=6cm]{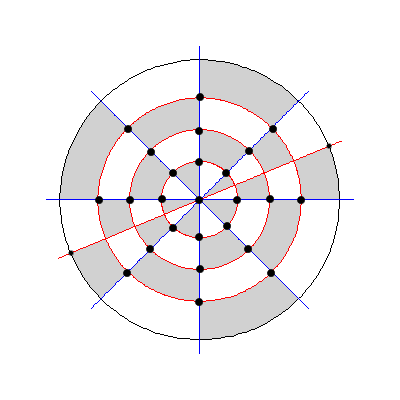}
\caption{Checkerboard: $\ell = 4$, $\mu > 0$}
\label{FEx2-nodal2}
\end{minipage} \hfill
\begin{minipage}[c]{.46\linewidth}
\includegraphics[width=6.5cm]{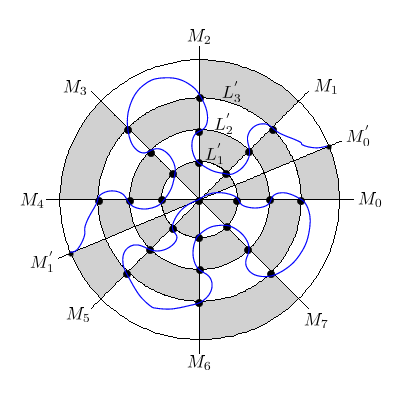}
\caption{$N(H^{\mu})$} \label{FEx2-nodalset}
\end{minipage}
\end{figure}

Figure~\ref{FEx2-nodalset} shows the nodal pattern of $H^{\mu}$ in
the exponential map, with one component tangent to the great circle
$M'_0 \cup M'_1$ at the north pole, the other at the south pole.

\section{Courant sharp property and open questions for minimal partitions for the sphere.}\label{s5}

Leydold's thesis \cite{Ley} (see also a preliminary analysis in
\cite{LeyD}) is devoted to this question. We reproduce below some
synthesis essentially extracted from \cite{HHOT1}. Given a spherical
harmonic $u$, let $\mu(u)$ denote the number of nodal domains of $u$
(this notation should not induce confusion with the parameter
$\mu$ appearing in the preceding sections). \vspace{-3mm}
\begin{itemize}
\item Courant's theorem for the sphere says that for any $u_{\ell}
\in \cH_{\ell}$,
$$
\mu(u_\ell) \leq \ell^2 +1\,,
$$
where the right-hand side is $1 + \sum_{k=0}^{\ell -1}\dim \cH_{k}$;
\item  Pleijel's asymptotic bound for the number of nodal domains
extends to bounded domains in $\mathbb{R}^n$, and more generally to
compact $n$-manifolds with boundary, with a universal constant $\gamma(n) <1$
replacing the constant $\gamma(2) = 4/(j_{0,1})^2$ in the right-hand
side of (\ref{Pl}) (Peetre \cite{Pe}, B\'{e}rard-Meyer \cite{BeMe}). It
is also interesting to note that this constant is independent of the
geometry. In particular Pleijel's theorem is true in the case of the
sphere. For any sequence of eigenfunctions $u_\ell \in \mathcal
H_\ell$
\begin{equation}\label{Pl}
\lim\sup_{\ell \rightarrow +\infty} \frac{\mu(u_\ell)}{\ell^2+1}
\leq \frac{4}{( j_{0,1})^2}\,.
\end{equation}
\item  Leydold stated the following conjecture on the maximal
cardinal of nodal sets of a spherical harmonic.
\begin{conjecture}\label{conjLey}
$$
\max_{u\in \mathcal H_\ell}\mu(u)= \left\{ \begin{array}{ll} \frac
12 (\ell +1)^2 & \mbox{ if } \ell \mbox{ is odd, }\\[4pt] \frac 12 \ell
(\ell +2) & \mbox{ if } \ell \mbox{ is even. }
\end{array}\right.
$$
\end{conjecture}

The values in the right hand side are the maximum numbers of
nodal domains of the decomposed spherical harmonics in spherical
coordinates. This  conjecture is proved in \cite{Ley} for
$\ell \leq 6$. Note that the example treated in Appendix~\ref{S-Maple} for
$\ell=3$ (middle subfigure in Fig. \ref{MEx1-L3-bifurc}) shows the
optimality in this case.   In \cite{LeyT}, Leydold constructs
regular spherical harmonics of degree $\ell$ with
$O(\frac{\ell^2}{4})$ nodal domains, see also
\cite[Theorem~2.1]{ErJaNa}.

Conjecture~\ref{conjLey} implies that the only Courant
sharp situations (that is situations in which Courant's upper bound is
attained in some eigenspace) correspond to the first and second
eigenvalues. This last statement is true as a consequence of
a theorem \`{a} la Courant, using the symmetry or antisymmetry of spherical
harmonics under the antipodal map, or as a corollary of the following theorem by
Karpushkin \cite{Ka}.

\begin{theorem}
$$
\max_{u\in \mathcal H_\ell}\mu (u) \le \left\{
\begin{array}{ll} \ell (\ell -2) +5  & \mbox{ if } \ell \mbox{ is
odd, }\\ \ell (\ell -2) +4 & \mbox{ if } \ell \mbox{ is even. }
\end{array}\right.
$$
\end{theorem}
Conjecture \ref{conjLey} implies the following inequality which
improves Pleijel's theorem.
\begin{conjecture}
For any sequence of eigenfunctions $u_\ell\in \mathcal H_\ell$, we have
\begin{equation}\label{conj5.3}
\lim\sup_{\ell \rightarrow +\infty} \frac{\mu(u_\ell)}{\ell^2+1} \leq  \frac 12\,.
\end{equation}
\end{conjecture}

It is easy to see that \eqref{conj5.3} cannot be improved  (look at
product eigenfunctions).

\item  Spectral minimal partitions are for example defined in \cite{HHOT}.
Motivated by a conjecture in harmonic analysis popularized by Bishop
\cite{Bis} (who refers to \cite{FrHa}), the authors of \cite{HHOT}
have proved in \cite{HHOT1}  that up to rotation the minimal
$3$-partition is the so-called $Y$-partition ($\{ 0< \phi
<\frac{2\pi}{3}\}$, $\{ \frac{2\pi}{3} < \phi <\frac{4\pi}{3}\}$,
and $\{ \frac{4\pi}{3} < \phi < 2\pi\}$). There is a conjecture that
the four faces of a spherical tetrahedron determine a  minimal
$4$-partition on $\mathbb S^2$. What we get from the previous item
and the general theory of \cite{HHOT}  (nodal minimal partitions
should correspond to a Courant sharp situation) is that minimal
$k$-partitions cannot be nodal for $k>2$.
\item With a different point of view, let us mention the contributions of
\cite{NaSo} on random spherical harmonics.
\end{itemize}

\appendix

\section{Some simulations with Maple}\label{S-Maple}

In this appendix, we provide some pictures issued from numerical
computations with Maple. The nodal sets are viewed in the
exponential map at the north pole. The outer circle, at distance
$\pi$, is the cut-locus of $p_{+}$ and corresponds to the south
pole.

\begin{figure}
\begin{minipage}[c]{.46\linewidth}
\includegraphics[width=7cm]{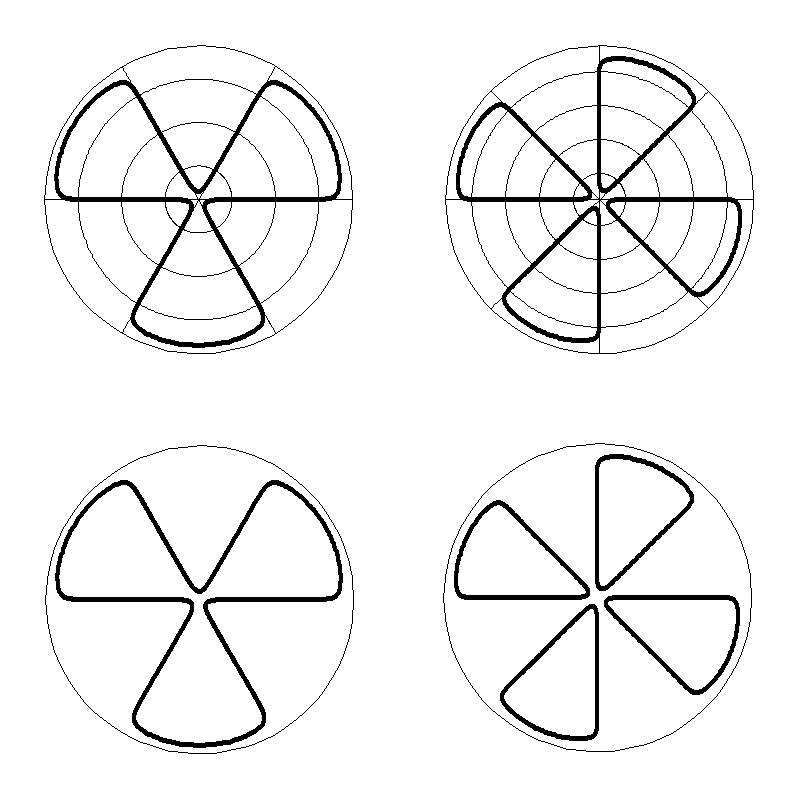}
\caption{Example 1: $\ell = 3$ and $\ell = 4$}\label{FMEx1-L3L4}
\end{minipage} \hfill
\begin{minipage}[c]{.46\linewidth}
\includegraphics[width=7cm]{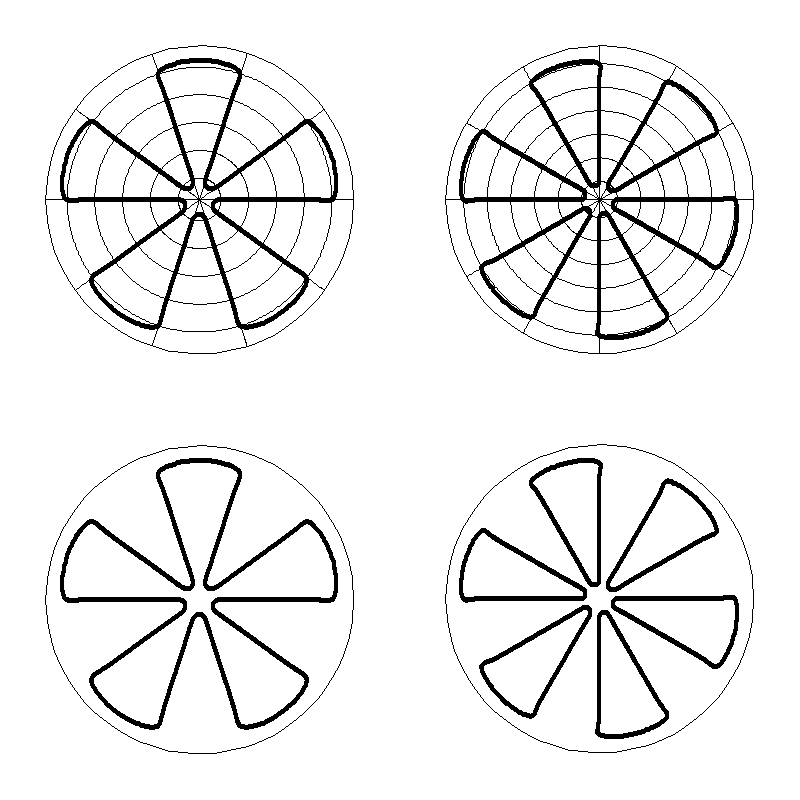}
\caption{Example 1: $\ell = 5$ and $\ell = 6$}\label{FMEx1-L5L6}
\end{minipage}
\end{figure}

Figures~\ref{FMEx1-L3L4} \resp \ref{FMEx1-L5L6} illustrate Proposition~\ref{Ex1-P} in the cases $\ell = 3$ (left) and $\ell = 4$ (right), \resp $\ell = 5$ (left) and $\ell = 6$ (right). They display the nodal set of $H^{\mu,\ell}$ (black thick line), with $\mu = 5\cdot 10^{-3}$ when $\ell$ is odd, and $\mu = 2\cdot 10^{-3}$ when $\ell$ even. The figures in the top line also display the checkerboards associated with $Z_{\ell}$ and $W_{\ell}$.

\begin{figure}
\centering
\includegraphics[width=0.8\linewidth]{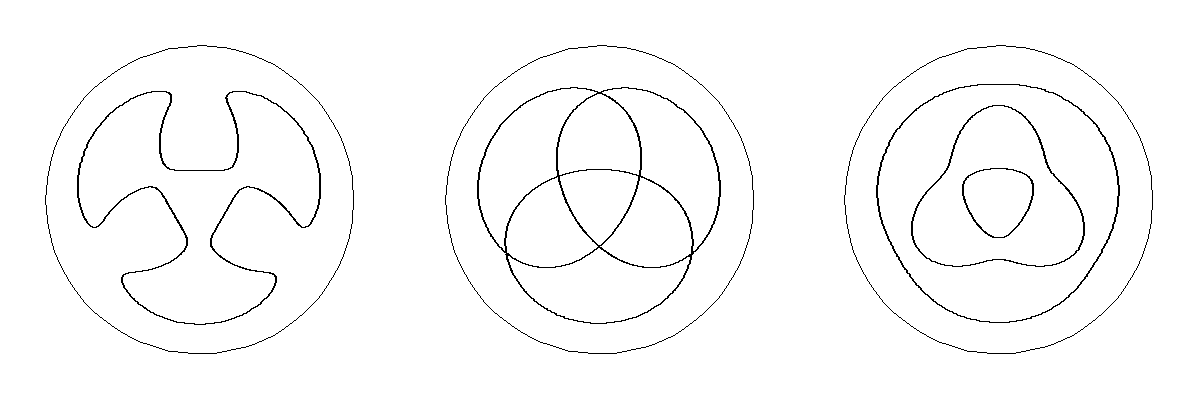}
\caption{Appearance and disappearance of critical zeros}\label{MEx1-L3-bifurc}
\end{figure}

Figure~\ref{MEx1-L3-bifurc} illustrates the occurrence of critical
zeros in Stern's Example 1, with $\ell = 3$. The corresponding
Legendre polynomial is $P_3(t) = \frac{1}{2} \, t \, (5 t^2-3)\,$. The
polynomial $P_2(t) = \frac{1}{2}\, (3t^2-1)$ has two roots $\pm \,
\frac{1}{\sqrt{3}}\,$. According to Section~\ref{SSS-Ex1-cz}, there
are twelve possible critical zeros, given in spherical coordinates
by the points $\left(\arccos(\pm
\frac{1}{\sqrt{3}}),j\frac{\pi}{6}\right)$ with $j \in \{1, 3, 5, 7,
9, 11\}$, and exactly two critical values of the parameter, $\mu =
\pm \,\sqrt{2}\,$. For $\mu > 0$, there is exactly one critical value
$\mu = \sqrt{2}$, which is associated with six critical zeros.
Figure~\ref{MEx1-L3-bifurc} shows the nodal set $N(H^{\mu,3})$ for
$\mu < \sqrt{2}$ (left), for $\mu = \sqrt{2}$ (center) and for $\mu
> \sqrt{2}$ (right).

\begin{figure}
\centering
\includegraphics[width=0.6\linewidth]{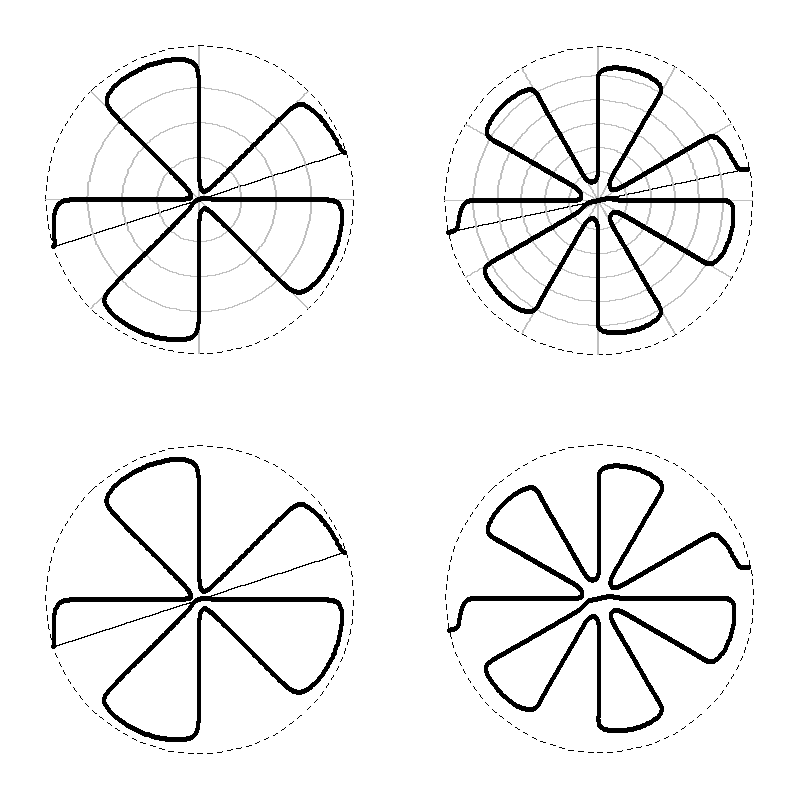}
\caption{Example 2: $\ell = 4$ and $\ell = 6$}\label{MEx2-L4L6}
\end{figure}

Figure~\ref{MEx2-L4L6} illustrates Proposition~\ref{Ex2-P}. The
figures display the nodal set of the function $H^{\mu}$ (thick lines), with $\mu = 10^{-3}$ and $\varepsilon = 0.4$. The figures in the top line also display the checkerboards associated with $W, V_{\alpha}$. The great circle $M'_0 \cup M'_1$ divides the sphere into two closed hemispheres. Each one contains a simple closed nodal curve tangent to the great circle at one of the poles. As usual, the south pole is represented by the outer circle (dotted line), the cut-locus of $0$ in the tangent space at $p_{+}$.

\newpage
\section{Translation of citations from Stern's thesis}\label{S-transl}

We provide below a rough translation of the citations from Stern's thesis in Section~\ref{S-intro}

\begin{quote}
[E1]~\ldots one can for example easily show that on the sphere the numbers $2$ or $3$ occur as number of nodal domains for each eigenvalue, and that for the square, if we arrange the eigenvalues in increasing order, the number $2$ always reappears as number of nodal domains.\medskip

[K1]~We shall then show that for each eigenvalue there exist eigenfunctions on the sphere whose nodal lines divide the sphere into   $2$ or $3$ nodal domains only \ldots The number of nodal domains $2$ occurs for the eigenvalues $\lambda_n = (2r+1)(2r+2)\, ~r=1,2,\cdots$; \\[3pt]
[K2]~similarly, we shall now show that $3$ occurs as number of nodal domains for all eigenvalues $ \lambda_n =2r  (2r+1), ~r=1,2,\cdots$.
\end{quote}

\begin{quote}
[I1]~Superimpose the systems of nodal lines of the two functions, and hatch the domains in which the functions have the same sign, then the nodal lines of the spherical function
$$
P^{2r+1}_{2r+1} (\cos \vartheta)\cos (2r+1) \varphi  + \mu
P_{2r+1}(\cos \vartheta)\,,\quad \mu >0
$$
can only pass in the non-hatched domains\\[3pt]
[I3]~ and hence for values of $\mu$ which are small enough, they stay in a neighborhood of the nodal lines of $$P^{2r+1}_{2r+1} (\cos
\vartheta)\cos (2r+1) \varphi,$$ \ie\, the $2r+1$ meridians, so that the system of nodal lines varies continuously with $\mu$ \ldots. \\[3pt]
[I2]~Furthermore, the nodal lines must pass through the $2(2r+1)^2$ intersection points of the system of nodal lines of the two above spherical functions \ldots
\end{quote}

\hfill \newpage

\section*{Acknowledgments} The authors would like to thank T.~Hoffmann-Ostenhof and J.~Leydold for useful discussions, and for providing the papers \cite{LeyD,LeyT}. They also thank the anonymous referee for his remarks.

\bibliographystyle{plain}

\begin{thebibliography}{}

\end{thebibliography}


\begin{thebibliography}{1}

\bibitem {BeMe} P. B\'erard and D. Meyer.
\newblock
In\'egalit\'es isop\'erim\'etriques et applications.
\newblock Annales scientifiques de l'\'Ecole Normale Sup\'erieure, 15:3 (1982), 513-541.


\bibitem{BeHe} P. B\'{e}rard and B. Helffer.
\newblock Dirichlet eigenfunctions of the square membrane:
Courant's property, and A.~Stern's  and  {\AA}.~Pleijel's analyses.
\newblock arXiv:14026054.
\newblock To appear in Springer Proceedings in Mathematics \& Statistics  (2015), MIMS-GGTM conference in memory of M.~S.~Baouendi. A. Baklouti, A. El Kacimi, S. Kallel, and N. Mir Editors.

\bibitem{Bis} C.J. Bishop.
\newblock Some questions concerning harmonic measure.
\newblock Dahlberg, B. (ed.) et al., Partial Differential equations
 with minimal smoothness and applications.
\newblock IMA Vol. Math. Appl. 42 (1992), 89-97.


\bibitem{Cou} R. Courant.
\newblock Ein allgemeiner Satz zur Theorie der Eigenfunktionen
selbstadjungierter Differentialausdr\"{u}cke.
\newblock Nachr. Ges. G\"{o}ttingen (1923), 81-84.

\bibitem{CH} R. Courant and D. Hilbert.
\newblock Methods of Mathematical Physics. Volume 1. John Wiley \&
Sons, 1989.

\bibitem{ErJaNa} A. Eremenko, D. Jakobson, and N. Nadirashvili.
\newblock On nodal sets and nodal domains on $\mathbb S^2$.
\newblock Annales Institut Fourier 57:7 (2007), 2345-2360.

\bibitem{FrHa} S. Friedland and W.K. Hayman.
\newblock Eigenvalue inequalities for the Dirichlet problem on spheres
 and the growth of subharmonic functions.
\newblock Comment. Math. Helvetici 51 (1976), 133-161.

 \bibitem{HHOT} B. Helffer, T. Hoffmann-Ostenhof and S. Terracini.
\newblock Nodal domains and spectral minimal partitions.
\newblock Annales Institut Henri Poincar\'e (Analyse non lin\'eaire)
26 (2009), 101-138.

\bibitem{HHOT1}
B. Helffer, T. Hoffmann-Ostenhof and S. Terracini.
\newblock On spectral minimal partitions: the case of the sphere.
\newblock Around the research of Vladimir Maz'ya. III,
p.~153-178, Int. Math. Ser. 13, Springer (N.Y.) 2010.

\bibitem{HP} B. Helffer and M. Persson Sundqvist.
\newblock Nodal domains in the square--the Neumann case--.
\newblock arXiv: 1410.6702.

\bibitem{Ka} V.N. Karpushkin.
\newblock Topology of the zeros of eigenfunctions.
\newblock Funktional Anal. i Prilozehen 23:3 (1989), 59-60.

\bibitem{Lew} H. Lewy.
\newblock On the minimum number of domains in which the nodal lines of
spherical harmonics divide the sphere.
\newblock Comm. Partial Differential Equations 2:12 (1977), 1233-1244.

\bibitem{LeyD} J. Leydold.
\newblock Knotenlinien und Knotengebiete von Eigenfunktionen.
\newblock Diplom Arbeit, Universit\"at Wien (1989), unpublished.
Available at \url{http://othes.univie.ac.at/34443/}


\bibitem{LeyT}
J. Leydold.
\newblock Nodal properties of spherical harmonics.
\newblock Dissertation Universit\"{a}t Wien (January 1993).


\bibitem{Ley}
J. Leydold.
\newblock On the number of nodal domains of spherical harmonics.
\newblock Topology 35 (1996), 301-321.


\bibitem{MOS} W.~Magnus, F.~Oberhettinger, and R.P.~Soni.
\newblock Formulas and Theorems for the Special Functions of Mathematical Physics.
Third Edition.
\newblock  Berlin: Springer-Verlag, 1966.

\bibitem{NaSo} F. Nazarov and  M. Sodin.
\newblock On the number of nodal domains of random spherical
harmonics.
\newblock Amer. J. Math. 131 (2009), 1337-1357.

\bibitem {Pe} J. Peetre.
\newblock A generalization of Courant nodal theorem.
\newblock Math. Scandinavica 5 (1957), 15-20.

\bibitem  {Pl} {\AA}.~Pleijel.
\newblock Remarks on Courant's nodal theorem.
\newblock Comm. Pure. Appl. Math. 9 (1956), 543-550.

\bibitem{Pol} I. Polterovich.
\newblock Pleijel's nodal domain theorem.
\newblock Proc. Amer. Math. Soc. 137 (2009), 1021-1024.


\bibitem{St} A. Stern.
\newblock Bemerkungen \"uber asymptotisches Verhalten von Eigenwerten
und Eigenfunktionen.
Inaugural-Dissertation zur Erlangung der Doktorw\"{u}rde der Hohen
Mathematisch-Naturwissenschaftlichen Fakult\"{a}t der Georg
August-Universit\"{a}t zu G\"{o}ttingen  (30 Juli 1924).
\newblock Druck der Dieterichschen Univertis\"{a}ts-Buchdruckerei (W. Fr. Kaestner).
G\"{o}ttingen, 1925

\bibitem{St1} A. Stern.
\newblock Bemerkungen \"uber asymptotisches Verhalten von Eigenwerten
und Eigenfunktionen.
\newblock Diss. G\"ottingen  (30 Juli 1924).
\newblock Extracts and annotations by P.~B\'{e}rard and B.~Helffer.
Available at
\url{http://www-fourier.ujf-grenoble.fr/~pberard/R/stern-1925-thesis-partial-reprod.pdf}

\bibitem{Sz} G. Szeg\"o.
\newblock Orthogonal Polynomials.
\newblock Fourth edition.  Amer. Math. Soc. Colloquium Publications,
Vol. XXIII, Amer. Math. Soc. Providence,
R.I. (1975).

\end{thebibliography}

\end{document}